\documentclass{degm}

\usepackage{amssymb}
\usepackage{amscd}
\usepackage[all]{xypic}

\newcommand{\Qp}{\mathbf{Q}_p}
\newcommand{\Zp}{\mathbf{Z}_p}
\newcommand{\Cp}{\mathbf{C}}
\newcommand{\Fp}{\mathbf{F}_p}
\newcommand{\eps}{\varepsilon}
\newcommand{\ra}{\rightarrow}

\newcommand{\Qpbar}{\overline{\mathbf{Q}}_p}
\newcommand{\Fpbar}{\overline{\mathbf{F}}_p}

\newcommand{\Kbar}{\overline{K}}
\newcommand{\KbarH}{\overline{K}^H}
\newcommand{\Fbar}{\overline{F}}
\newcommand{\bnrigplus}[1]{\mathbf{B}^+_{\mathrm{rig} #1}}

\newcommand{\on}{\operatorname}
\newcommand{\OO}{\mathcal{O}}

\renewcommand{\hat}{\widehat}
\renewcommand{\tilde}{\widetilde}
\renewcommand{\phi}{\varphi}
\renewcommand{\projlim}{\varprojlim}
\newcommand{\ZZ}{\mathbf{Z}}

\newcommand{\QQ}{\mathbf{Q}}
\newcommand{\FF}{\mathbf{F}}
\newcommand{\BB}{\mathbf{B}}

\newcommand{\btst}[2]{\widetilde{\mathbf{B}}^{\dagger #1}_{\mathrm{log} #2}}
\newcommand{\btrig}[2]{\widetilde{\mathbf{B}}^{\dagger #1}_{\mathrm{rig} #2}}
\newcommand{\bnst}[2]{\mathbf{B}^{\dagger #1}_{\mathrm{log} #2}}
\newcommand{\bcris}{\mathbf{B}_{\mathrm{cris}}}
\newcommand{\bnrig}[2]{\mathbf{B}^{\dagger #1}_{\mathrm{rig} #2}}
\newcommand{\btstplus}[1]{\widetilde{\mathbf{B}}^{+}_{\mathrm{log} #1}}
\newcommand{\btrigplus}[1]{\widetilde{\mathbf{B}}^{+}_{\mathrm{rig} #1}}

\newcommand{\dnst}[1]{\mathbf{D}^{\dagger #1}_{\mathrm{log}}}
\newcommand{\dnrig}[1]{\mathbf{D}^{\dagger #1}_{\mathrm{rig}}}

\newcommand{\bmax}{\mathbf{B}_{\mathrm{max}}}

\newcommand{\bst}{\mathbf{B}_{\mathrm{st}}}
\newcommand{\bdR}{\mathbf{B}_{\mathrm{dR}}}
\newcommand{\btdag}[1]{\widetilde{\mathbf{B}}^{\dagger #1}}
\newcommand{\bdag}[1]{\mathbf{B}^{\dagger #1}}
\newcommand{\bplus}{\mathbf{B}^+}
\newcommand{\btplus}{\widetilde{\mathbf{B}}^+}
\newcommand{\bt}{\widetilde{\mathbf{B}}}

\newcommand{\atdag}[1]{\widetilde{\mathbf{A}}^{\dagger #1}}
\newcommand{\adag}[1]{\mathbf{A}^{\dagger #1}}
\newcommand{\aplus}{\mathbf{A}^+}
\newcommand{\atplus}{\widetilde{\mathbf{A}}^+}
\newcommand{\at}{\widetilde{\mathbf{A}}}
\newcommand{\bhol}[1]{\mathbf{B}^+_{\mathrm{rig} #1}}
\newcommand{\blog}[1]{\mathbf{B}^+_{\mathrm{log} #1}}
\newcommand{\e}{\mathbf{E}}

\newcommand{\etplus}{\widetilde{\mathbf{E}}^+}
\newcommand{\et}{\widetilde{\mathbf{E}}}
\newcommand{\dst}{\mathbf{D}_{\mathrm{st}}}
\newcommand{\dcris}{\mathbf{D}_{\mathrm{cris}}}
\newcommand{\ddR}{\mathbf{D}_{\mathrm{dR}}}
\newcommand{\dfont}{\mathbf{D}}
\newcommand{\dfontp}{\mathbf{D'}}
\newcommand{\ddif}{\mathbf{D}_{\mathrm{dif}}}
\renewcommand{\ddag}[1]{\mathbf{D}^{\dagger #1}}
\newcommand{\dsen}{\mathbf{D}_{\mathrm{Sen}}}
\newcommand{\bsen}{\mathbf{B}_{\mathrm{Sen}}}

\newcommand{\vale}{v_\mathbf{E}}

\newcommand{\ndr}{\mathbf{N}_{\mathrm{dR}}}

\newenvironment{bibli}{}

\title{An introduction to the theory of $p$-adic representations}

\author{Laurent Berger}

\date{October 2002}

\contact[\texttt{laurent@math.harvard.edu}]{Laurent Berger \\
          Harvard University \\
          Department of Mathematics \\
          One Oxford Street \\
          Cambridge, MA 02138-2901 \\ USA}

\begin{document}

\maketitle

\newcommand{\Subsection}{\subsection}
\renewcommand{\thesection}{\Roman{section}}
\renewcommand{\thesubsection}{\Roman{section}.\arabic{subsection}}

\setcounter{tocdepth}{3}

\begin{abstract}
The purpose of this informal article is to 
introduce the reader to some of the objects 
and methods of the theory of $p$-adic representations. 
My hope is that students and
mathematicians who are new to the subject
will find it useful as a starting point.
It consists mostly of an expanded version of
the notes for my two lectures at the 
``Dwork trimester'' in June 2001.
\end{abstract}


\tableofcontents

\[ \spadesuit\heartsuit\clubsuit\diamondsuit \]

\newpage

\section{Introduction}
\Subsection{Introduction}
\subsubsection{Motivation}\label{motiv}
One of the aims of arithmetic geometry is to understand the structure of the Galois group
$\on{Gal}(\overline{\QQ}/\QQ)$, or at least to understand its action on representations
coming from geometry. A good example is the Tate module $T_{\ell}E$ of an elliptic curve
$E$ defined over $\QQ$. The action of $\on{Gal}(\overline{\QQ}/\QQ)$ on $T_{\ell}E$ carries a
lot of arithmetical information, including the
nature of the reduction of $E$ at various primes and the number of
points in $E(\FF_q)$.

Let  $D_p \subset \on{Gal}(\overline{\QQ}/\QQ)$ be 
the decomposition group of a place above $p$; it is naturally isomorphic to
$\on{Gal}(\Qpbar/\Qp)$.
The aim of the theory of $p$-adic representations is to extract
information from the action of $D_p$, on
$\Qp$-vector spaces. This is in stark contrast to 
the theory of $\ell$-adic representations, which
endeavors to understand the action of $D_p$ on $\QQ_{\ell}$-vector spaces with $\ell\neq p$.

In this latter
situation, the topology of $D_p$ is mostly incompatible with that of an $\ell$-adic vector
space (essentially because the wild inertia is a $p$-group), 
and the result is that the theory of $\ell$-adic representations is of an algebraic
nature. On the other hand, in the $p$-adic case, the topologies are compatible and as a
result there are far too many representations. The first step is therefore to single out the
interesting objects, and to come up with significant invariants attached to them. Unlike the
$\ell$-adic situation, the study of $p$-adic representations is therefore of a rather
($p$-adic) analytic nature. 

For example, there exists a $p$-adically
continuous family of characters
of the group 
$\on{Gal}(\Qpbar/\Qp)$, given by $\chi^s$ where $\chi$ is the cyclotomic character and
$s$ varies in weight space (essentially $p-1$ copies of $\Zp$). Out of those characters,
only those corresponding to integer values of $s$ ``come from geometry''. This kind of
phenomenon does not arise in the $\ell$-adic case, where every character is ``good''.

The aim of this article is to introduce some of the objects and techniques which are 
used to study $p$-adic representations, and to provide explanations of recent developments.

\subsubsection{Organization of the article} 
This article is subdivided in chapters, each of which is subdivided 
in sections made up of paragraphs. At the
end of most paragraphs, I have added references to the literature. This article is as
much a quick survey of some topics as a point of entry for the literature on the subject.

\begin{bibli}
References are indicated at the end of paragraphs. For each topic, I have tried to
indicate a sufficient number of places where the reader can find all
the necessary details.
I have not always tried to give references to original articles, but
rather to more recent (and sometimes more readable) accounts. 
\end{bibli}

\subsubsection{Acknowledgments} The basis for this article are the
two lectures which I gave at the ``Dwork Trimester'' in Padova, and I
thank the organizers, especially F. Baldassarri, 
P. Berthelot and B. Chiarellotto
for the time and effort they spent to make this
conference a success. After I wrote a first version of this article,
M. \c{C}iperjani, J-M. Fontaine and H. J. Zhu took the time to read it,
pointed out several inaccuracies and
made many suggestions for improvement. 
Any remaining inaccuracies are
entirely my fault. 

\Subsection{$p$-adic representations}
\subsubsection{Some notations}\label{nota}
The results described in this article are true in a rather general setting. Let $k$ be
any perfect field of characteristic $p$ 
(\emph{perfect} means that the map $x \mapsto x^p$ is an
automorphism), and let $F=W(k)[1/p]$ be the fraction field of $\OO_F=W(k)$, the ring of Witt
vectors over $k$
(for reminders on Witt vectors, see paragraph \ref{witt}).
Let $K$ be a totally ramified extension of $F$, and let $\Cp=\hat{\Fbar}=\hat{\overline{K}}$
be the $p$-adic completion of the algebraic closure of $F$. If $k$ is
contained in the algebraic closure of $\Fp$, then $\Cp=\Cp_p$, the
field of so-called $p$-adic complex numbers. 

An important special case is when $k$ is a finite extension of $\Fp$, so that 
$K$ is a finite extension of $\Qp$, and $F$ is then 
the maximal unramified extension of $\Qp$
contained in $K$. The reader can safely assume that we're 
in this situation throughout the
article. Another important special case though
is when $k$ is algebraically closed.

We'll choose once and for all a compatible sequence of primitive $p^n$-th roots of unity,
$\eps^{(0)}=1$, and $\eps^{(n)} \in \mu_{p^n} \subset \Kbar$, 
such that $\eps^{(1)} \neq 1$ and
$(\eps^{(n+1)})^p =
\eps^{(n)}$. Let $K_n=K(\eps^{(n)})$ and $K_{\infty} = \cup_{n=0}^{+\infty} K_n$. Making
such a choice of $\eps^{(n)}$ is like choosing an orientation 
in $p$-adic Hodge theory, in the
same way that choosing one of $\pm i$ is 
like choosing an orientation in classical geometry. Here are the
various fields that we are considering:
\[ F \subset \overset{G_K}
{\overbrace{\underset{\Gamma_K}
{\underbrace{K \subset K_n \subset K_\infty}} = 
\underset{H_K}{\underbrace{
  K_\infty \subset \overline{F} = \overline{K}}}}}
\subset \Cp \]

Let $G_K$ be the Galois group $\on{Gal}(\overline{K}/K)$. The \emph{cyclotomic 
character} $\chi:G_K \ra \Zp^*$ is defined by
$\sigma(\zeta)=\zeta^{\chi(\sigma)}$ for every $\sigma \in G_K$ and
$\zeta \in \mu_{p^{\infty}}$.
The kernel of the cyclotomic character is $H_K=\on{Gal}(\overline{K}/K_{\infty})$, and
$\chi$ therefore identifies 
$\Gamma_K=\on{Gal}(K_{\infty}/K)=G_K/H_K$ with an open subgroup of $\Zp^*$.

\subsubsection{Definitions}\label{defi}
A \emph{$p$-adic representation} $V$ is a finite dimensional $\Qp$-vector space with a
continuous linear action of $G_K$. 
The dimension of $V$ as a $\Qp$-vector space will always be denoted by $d$.
Here are some examples of $p$-adic representations:
\begin{enumerate}
\item If $r\in\ZZ$, then $\Qp(r)=\Qp\cdot e_r$ where $G_K$ acts on $e_r$ by
$\sigma(e_r)=\chi(\sigma)^r e_r$. This is the $r$-th Tate twist of $\Qp$;
\item if $E$ is an elliptic curve, then the Tate module of $E$, 
$V=\Qp \otimes_{\Zp} T_pE$ is a
$p$-adic representation of dimension $d=2$;
\item more generally, if $X$ is a proper and smooth variety over $K$, then the \'etale
cohomology $H^i_{\mathrm{\acute{e}t}}(X_{\overline{K}},\Qp)$ is a $p$-adic representation.
\end{enumerate}

This last example is really the most interesting
(the first two being special cases),
and it was the motivation for the
systematic study of $p$-adic representations. 
Grothendieck had suggested in 1970 the existence of a ``mysterious
functor'' (\emph{le foncteur myst\'erieux})
linking the \'etale and crystalline cohomologies of a
$p$-divisible group. Fontaine gave an algebraic construction 
of that functor which conjecturally allowed one to
recover,  for any $i$ and any proper and smooth $X/K$, 
the de Rham cohomology of $X/K$ (which is a filtered $K$-vector space) 
from the data of $H^i_{\mathrm{\acute{e}t}}(X_{\overline{K}},\Qp)$ 
as a $p$-adic representation.
His construction was shown to be valid by Tsuji; we'll come back to
that in \ref{tsuji}.

The above result is a $p$-adic analogue of the rather easy 
fact that if $X$ is a proper smooth variety over a number
field $L$, then over the complex numbers $\mathbb{C}$
one has an isomorphism
\[ \mathbb{C} \otimes_L H^i_{dR}(X/L) \simeq \mathbb{C} \otimes_{\ZZ} H^i(X,\ZZ) \]
given by integrating differential forms on cycles.

\subsubsection{Fontaine's strategy}\label{strat}
Fontaine's strategy for studying $p$-adic representations was to
construct \emph{rings of periods}, which are
topological
$\Qp$-algebras $B$, with an action of $G_K$ and some additional structures which 
are compatible with the action of
$G_K$ (for example: a Frobenius $\phi$, a filtration $\on{Fil}$, a
monodromy map $N$, a differential operator $\partial$), such that the $B^{G_K}$-module
$D_B(V)=(B \otimes_{\Qp} V)^{G_K}$, which inherits the additional
structures, is an interesting invariant of $V$.  
For Fontaine's constructions to work, one needs to assume that $B$ is
\emph{$G_K$-regular}, meaning that if $b \in B$ is such that the line $\Qp \cdot b$ is
stable by $G_K$, then $b \in B^*$. In particular, $B^{G_K}$ has to be a field.

In general, 
a simple computation shows that $\dim_{B^{G_K}} D_B(V) \leq d = \dim_{\Qp} V$,
and we say that $V$ is 
\emph{$B$-admissible} if equality holds, which is equivalent to having
$B\otimes_{\Qp}V \simeq B^d$ as $B[G_K]$-modules.
In this case, $B \otimes_{B^{G_K}} D_B(V) \simeq B \otimes_{\Qp} V$, and the
coefficients of a matrix of this isomorphism in two basis of $D_B(V)$ and $V$ are called
the \emph{periods} of $V$. 

Let us briefly mention a cohomological version of this: a $p$-adic representation $V$
determines a class $[V]$ in $H^1(G_K,\on{GL}(d,\Qp))$, and therefore a class $[V]_B$ in
$H^1(G_K,\on{GL}(d,B))$. The representation $V$ is $B$-admissible if
and only if $[V]_B$ is trivial.
In this case, $[V]_B$ is a coboundary, given explicitly by writing down a $G_K$-invariant
basis of $B \otimes_{\Qp} V$.

Here are some examples of rings of periods:
\begin{enumerate}
\item If $B=\overline{K}$, then $B^{G_K}=K$ and $V$ is $\overline{K}$-admissible if and only
if the action of $G_K$ on $V$ factors through a finite quotient. This is essentially Hilbert
90;
\item If $B=\Cp$, then $B^{G_K}=K$ (the so-called theorem of
  Ax-Sen-Tate, first shown by Tate)
and $V$ is $\Cp$-admissible if and only
if the action of the inertia $I_K$ on $V$ factors through a finite quotient. This was
conjectured by Serre and proved by Sen. We will return to this in \ref{tsuzuki};
\item Let $B=\bdR$ be Fontaine's ring of $p$-adic periods 
(defined below in \ref{bdr}). 
It is a field, equipped with a
filtration, and  $B^{G_K}=K$. If
$V=H^i_{\mathrm{\acute{e}t}}(X_{\overline{K}},\Qp)$, 
for a proper smooth $X/K$, then
$V$ is $\bdR$-admissible, and $\ddR(V) 
=\mathbf{D}_{\bdR}(V)
\simeq H^i_{\mathrm{dR}}(X/K)$ as filtered
$K$-vector spaces. This is one of the most important theorems of $p$-adic Hodge theory.
\end{enumerate}

\begin{bibli}
For the yoga of rings of periods and Tannakian categories, see Fontaine's \cite{Bu88sst}.
\end{bibli}

\Subsection{Fontaine's classification}
By constructing many rings of periods,
Fontaine has defined several subcategories of the category of all $p$-adic representations,
and in this paragraph, we shall list a number of them along with categories of invariants
which one can attach to them. 
Many of the words used here will be defined later in the text,
but the table below should serve as a guide to the world of $p$-adic representations.

{\footnotesize \begin{center}
\begin{tabular}{||p{18mm}|p{50mm}|p{16mm}|p{22mm}|} 
\hline
$p$-adic representations & Invariants & References in the text & $\ell$-adic analogue \\
\hline \hline
all of them  & $(\phi,\Gamma)$-modules & \ref{pgmod} &
 -- \\
\hline
Hodge-Tate & Hodge-Tate weights & \ref{sen} & -- \\
\hline
de Rham & 1. $p$-adic differential equations 

2. filtered $K$-vector spaces & \ref{bdr},  \ref{ndr} &
all $\ell$-adic \\
\hline
potentially semi-stable & quasi-unipotent differential equations 
& \ref{bst},  \ref{diffst} & quasi-unipotent \\
\hline
semi-stable  & 1. unipotent differential equations

2. admissible
filtered $(\phi,N)$-modules & \ref{bst},  \ref{diffst} &
unipotent \\
\hline
crystalline &  1. trivial differential equations

2. admissible filtered $\phi$-modules & \ref{bcris}, \ref{diffst} &
good reduction \\
\hline
\end{tabular}
\end{center}}

Each category of representations is a subcategory of the one above it.
One can associate to every $p$-adic representation a $(\phi,\Gamma)$-module, which is an
object defined on the boundary of the open unit disk. This object extends to a small
annulus, and if $V$ is de Rham, the action of the Lie algebra of $\Gamma$ gives a $p$-adic
differential equation. This equation is unipotent exactly when $V$ (restricted to
$G_{K_n}$ for some $n$) is semi-stable. 
In this case, the kernel of the connexion is a
$(\phi,N)$-module which coincides with the $(\phi,N)$-module attached to $V$ by $p$-adic
Hodge theory (one loses the filtration, however).

All of this will be explained later in the body of the text.

\[ \diamondsuit\spadesuit\heartsuit\clubsuit \]

\newpage

\section{$p$-adic Hodge theory}
In this chapter, we'll define the various rings of periods which are used in $p$-adic Hodge
theory, and give some simple examples of Fontaine's construction for an explicit geometric
object (an elliptic curve).

\Subsection{The field $\Cp$ and the theory of Sen}
Before we define the rings of periods which are used in $p$-adic Hodge theory, we'll review
some simple properties of the field $\Cp$ of $p$-adic complex numbers. As we saw above,
$\Cp$ is not a great ring of periods (since $\Cp$-admissible representations are potentially
unramified - and representations coming from arithmetic 
geometry are much more complicated than that), 
but one can still extract a lot of arithmetic
information from the data of $\Cp \otimes_{\Qp} V$: this is the
content of \emph{Sen's theory}.

\subsubsection{The action of $G_K$ on $\Cp$}\label{ast}
An important property of $\Cp$ that we will need is the computation of $\Cp^H$ where $H$
is a closed subgroup of $G_K$. Clearly, $\KbarH \subset \Cp^H$ and therefore
$\hat{\KbarH} \subset \Cp^H$. The Ax-Sen-Tate theorem says that the latter inclusion is
actually an equality: $\hat{\KbarH} = \Cp^H$. This was first shown by
Tate, and the proof was later improved and generalized by Sen and Ax.
Following Sen, Ax gave a natural proof of that result,
by showing that an element of $\Kbar$ almost invariant by $H$ is almost in $\KbarH$. 

The first indication that $\Cp$ was not a good choice for a 
ring of periods was given by a theorem of
Tate, which asserts that $\Cp$ does not contain periods for characters which are
too ramified (for example: 
the cyclotomic character). More precisely, he showed that if $\psi :
G_K \ra \Zp^*$ is a character which is trivial on $H_K$ but
which does not factor through a finite quotient of $\Gamma_K$, then
\[ H^0(K,\Cp(\psi^{-1})) = \{ x \in \Cp,\ g(x) = \psi(g)x\ \text{for
  all}\ 
g\in G_K\}=\{0\}. \] 
In particular, there is no
period in $\Cp$ for the cyclotomic character 
(a non-zero element of the above set is a period for $\psi^{-1}$).
Let us explain the proof of Tate's result; 
by the Ax-Sen-Tate theorem, the invariants of $\Cp$ under the action of $H_K$ are given by
$\Cp^{H_K}=\hat{K}_{\infty}$. The main argument in Tate's proof is the construction of
generalized trace maps $\on{pr}_{K_n} : 
\hat{K}_{\infty} \ra K_n$. The map $\on{pr}_{K_n}$ is a
continuous, $K_n$-linear, and $\Gamma_K$-equivariant
section of the inclusion $K_n \subset \hat{K}_{\infty}$.
In addition, if $x \in \hat{K}_{\infty}$, then $x = \lim_{n \ra \infty}
\on{pr}_{K_n}(x)$. 
We see that we can and should set $\on{pr}_{K_n}(x) = \lim_{m \ra +\infty}
[K_{n+m}:K_n]^{-1} \on{Tr}_{K_{n+m}/K_n} (x)$. The proof of the 
convergence of the above limit depends essentially on a good understanding of
the ramification of $K_{\infty}/K$.

Using these maps, one can prove Tate's theorem. Let $x$ be a period of
$\psi$. Since $\psi_{|H_K} = 1$, one has $x \in \Cp^{H_K}=\hat{K}_{\infty}$. We therefore
have $x =  \lim_{n \ra \infty} x_n$ 
where $x_n  = \on{pr}_{K_n}(x)$, and since $g(x)-\psi(g)x=0$ for
all $g \in G_K$, and $\on{pr}_{K_n}$ is Galois-equivariant, one also has 
$g(x_n)-\psi(g)x_n=0$ for all $g \in G_K$. If $x_n \neq 0$, this would imply that $\psi$
factors through $\on{Gal}(K_n/K)$, a contradiction, 
so that $x_n=0$ for every $n$. Since $x=\lim_{n \ra
\infty} x_n$, we also have $x=0$.

\begin{bibli}
General facts on $\Cp$ can be found in Koblitz's \cite{NK84}, which is a good introduction 
to $p$-adic numbers. The beginning of \cite{DGS} is a wonderful
introduction too.
The proof of Ax-Sen-Tate's theorem that we referred to
is in Ax's \cite{Ax70}, see also Colmez's \cite[\S 4]{Co00}. Tate's
theorems on the cohomology of $\Cp$ are in \cite{Ta62} or 
in Fontaine's \cite[\S 1]{Fo00}.
\end{bibli}

\subsubsection{Sen's theory}\label{sen}
The point of Sen's theory is to study the residual action of $\Gamma_K$ on the
$\hat{K}_{\infty}$-vector space $(\Cp \otimes_{\Qp} V)^{H_K}$, where
$V$ is any $p$-adic representation of $G_K$. 
We can summarize his main
result as follows. If $d \geq 1$, then 
$H^1(H_K,\on{GL}(d,\Cp))$ is trivial and 
the natural map: $H^1(\Gamma_K,\on{GL}(d,K_{\infty})) \ra
H^1(\Gamma_K,\on{GL}(d,\hat{K}_{\infty}))$
induced by the inclusion $K_{\infty} \subset \hat{K}_{\infty}$ is a bijection.

One can show that this implies the following:
given a $p$-adic representation $V$, the 
$\hat{K}_{\infty}$-vector space $(\Cp \otimes_{\Qp} V)^{H_K}$ has dimension
$d=\dim_{\Qp}(V)$, and the union of the finite dimensional $K_{\infty}$-subspaces of 
$(\Cp \otimes_{\Qp} V)^{H_K}$ stable by the action of $\Gamma_K$ 
is a $K_{\infty}$-vector space of dimension $d$. We shall
call it $\dsen(V)$, and the natural map $\hat{K}_{\infty} \otimes_{K_{\infty}} \dsen(V) 
\ra (\Cp \otimes_{\Qp} V)^{H_K}$ is
then an isomorphism. The $K_{\infty}$-vector space $\dsen(V)$ is endowed with an action
of $\Gamma_K$, and 
\emph{Sen's invariant} is the linear map giving the action of
$\on{Lie}(\Gamma_K)$ on $\dsen(V)$. It is the operator defined 
in $\on{End}(\dsen(V))$ by $\Theta_V =
\log(\gamma)/\log_p(\chi(\gamma))$, where $\gamma \in \Gamma_K$ is close enough to $1$ (the
definition of $\Theta_V$ obviously doesn't depend upon the choice of $\gamma$).

The operator $\Theta_V$ is then an invariant
canonically attached to $V$. Let us give a few examples: we say that
$V$ is \emph{Hodge-Tate}, with
\emph{Hodge-Tate weights} $h_1,\cdots,h_d \in \ZZ$, 
if there is a decomposition of $\Cp[G_K]$-modules: 
$\Cp \otimes_{\Qp} V = \oplus_{j=1}^d \Cp(h_j)$. Clearly, this is equivalent to $\Theta_V$
being diagonalizable with integer eigenvalues. For this reason, the eigenvalues of
$\Theta_V$ are usually called the 
\emph{generalized Hodge-Tate weights} of $V$. All representations
coming from a proper smooth 
variety $X/K$ (the subquotients of its \'etale cohomology groups)
are Hodge-Tate, and the $h_j$'s are the opposites of the jumps of
the filtration on the de Rham cohomology of $X$. For example, the Hodge-Tate weights of
$V=\Qp \otimes_{\Zp} T_p E$, where $E$ is an elliptic curve, are $0$ and $1$. 
Here is a representation which is not Hodge-Tate: let $V$ be a two dimensional
$\Qp$-vector space on which $G_K$ acts by
\[ \begin{pmatrix}
1 & \log_p(\chi(\overline{g})) \\
0 & 1
\end{pmatrix}
\qquad\text{so that}\qquad
\Theta_V = 
\begin{pmatrix}
0 & 1 \\
0 & 0
\end{pmatrix}. \]

\begin{bibli}
Relevant papers of Sen are \cite{Sn72,Sn73,Sn80}
and \cite{Sn93} which deals with families of representations
(about that, see also \ref{famil}). 
Colmez has given a different construction
more in the spirit of the ``ring of periods'' approach (by constructing a ring $\bsen$), 
see \cite{Co94}. For an interesting discussion of
all this, see Fontaine's course \cite[\S 2]{Fo00}.
\end{bibli}

\Subsection{The field $\bdR$}
\subsubsection{Reminder: Witt vectors}\label{witt}
Before we go any further, we'll briefly review the theory of Witt vectors. Let $R$ be a
perfect ring of characteristic $p$. For example, $R$ could be a finite field or an
algebraically closed field, or the ring of integers of an algebraically closed field (in
characteristic $p$, of course). The aim of the theory of Witt vectors is to construct a ring
$A$, in which $p$ is not nilpotent, and such that $A$ is separated and complete for the
topology defined by the ideals $p^n A$. We say that $A$ is a 
\emph{strict $p$-ring} with residual
ring $R$. The main result is that if $R$ is a perfect ring of characteristic $p$, then there
exists a unique (up to unique isomorphism) strict $p$-ring $A=W(R)$
with residual ring $R$.
It is called the ring of \emph{Witt vectors} over $A$.
Furthermore, because of the unicity, if one has a map $\xi: R \ra S$, then it lifts to a map
$\xi: W(R) \ra W(S)$. In particular, the map $x \ra x^p$ lifts to a Frobenius automorphism
$\phi: W(R) \ra W(R)$. 

Let us give a few simple examples: if $R=\Fp$, then $W(R)=\Zp$ and more generally, if $R$ is
a finite field, then $W(R)$ is the ring of integers of the unique unramified extension of
$\Qp$ whose residue field is $R$. If $R=\Fpbar$, then
$W(R)=\OO_{\hat{\Qp^{\mathrm{unr}}}}$. 
In the following paragraphs, we will see more interesting examples.

If $x=x_0 \in R$, 
then for every $n\geq 0$, choose an element
$\tilde{x}_n$ in $A$ whose image
in $R$ is $x^{p^{-n}}$. The sequence $\tilde{x}_n^{p^n}$ then converges in $A$ to an element
$[x]$ which depends only on $x$. This defines a multiplicative map $x \mapsto [x]$ from $R
\ra A$, which is a section of the projection
$x \mapsto \overline{x}$, called the \emph{Teichm\"uller map}. These
Teichm\"uller elements are a distinguished set of representatives of the elements of $R$.
One can write every element $x \in A$ in a unique way as $x=\sum_{n=0}^{+\infty} p^n [x_n]$
with $x_n \in R$.
Given two elements $x, y \in A$, one can then write 
\[ x+y  = \sum_{n=0}^{+\infty} p^n [S_n(x_i,y_i)] \quad\text{and}\quad
xy  = \sum_{n=0}^{+\infty} p^n [P_n(x_i,y_i)] \]
where $S_n$ and $P_n \in \ZZ[X_i^{p^{-n}},Y_i^{p^{-n}}]_{i=0 \cdots n}$
are universal homogeneous polynomials of degree one
(if one decides that the degrees of the $X_i$ and $Y_i$ are $1$).
For example, $S_0(X_0,Y_0)=X_0+Y_0$ and
$S_1(X_0,X_1,Y_0,Y_1)=X_1+Y_1+p^{-1}((X_0^{1/p}+Y_0^{1/p})^p-X_0-Y_0)$. 
The simplest way to construct $W(R)$ is then by setting
$W(R)=\prod_{n=0}^{+\infty} R$ and by making it into a ring using the addition and
multiplication defined by the $P_n$ and $S_n$, which are given by 
(not so) simple 
functional equations.

Finally, let us mention that if $R$ is not perfect, then there still exist strict $p$-rings
$A$ such that $A/pA=R$, but $A$ is not unique anymore. Such a ring is
called a 
\emph{Cohen ring}.
For example, if $R=\Fp[[X]]$, then one can take $A=\Zp[[X]]$, but for all $\alpha \in p\Zp$,
the map $X \mapsto X + \alpha$ is a non-trivial isomorphism of $A$ which induces the
identity on $R$.

\begin{bibli}
The above summary is inspired from a course by P. Colmez. The best place to 
start further reading is Harder's survey \cite{GH97}. The construction of Witt vectors is
also explained by Serre in \cite{SeCL} (or in English in \cite{SeCLeng}).
\end{bibli}

\subsubsection{The universal cover of $\Cp$}\label{thickcp}
Let $\etplus$ be the set defined by
\[ \etplus=\projlim_{x\mapsto x^p} \OO_{\Cp} 
=\{ (x^{(0)},x^{(1)},\cdots) \mid (x^{(i+1)})^p = x^{(i)} \} \]
which we make into a ring by deciding that if $x=(x^{(i)})$ and $y=(y^{(i)})$ 
are two elements of $\etplus$, then their sum and their product is defined by:
\[ (x+y)^{(i)}= \lim_{j \ra \infty} (x^{(i+j)}+y^{(i+j)})^{p^j} \ \text{and}\ 
(xy)^{(i)}=x^{(i)}y^{(i)}. \] 
This makes $\etplus$ into a characteristic $p$ local 
ring. Let $\eps=(\eps^{(i)})$ where the $\eps^{(i)}$ are the elements which 
have been chosen in \ref{nota}. It is easy to see that 
$\Fp((\eps-1)) \subset \et = \etplus[(\eps-1)^{-1}]$ and one can show that 
$\et$ is a field which is the completion of
the algebraic (non-separable!) closure of $\Fp((\eps-1))$, so it is
really a familiar object. 

We define a 
valuation $\vale$ on $\et$ by $\vale(x)=v_p(x^{(0)})$ so that
$\etplus$ is the integer ring of $\et$ for $\vale$. For example,
$\vale(\eps-1)= \lim_{n \ra \infty} v_p (\eps^{(n)}-1)^{p^n}=p/(p-1)$.

Finally, let us point out that 
there is a natural map $\etplus \ra \projlim_{x\mapsto x^p} \OO_{\Cp}/p$ and
it's not hard to show that this map is an isomorphism.

There is a natural map $\theta$ from $\etplus$ to $\OO_{\Cp}$, which sends $x=(x^{(i)})$
to $x^{(0)}$, and the map $\theta: \etplus \ra \OO_{\Cp}/p$ is a homomorphism.
Let $\atplus=W(\etplus)$ and \[ \btplus =\atplus[1/p] 
=\{ \sum_{k\gg -\infty} p^k [x_k],\ x_k \in \etplus
\}\] where $[x] \in \atplus$ is the Teichm{\"u}ller lift of $x \in \etplus$.
The map $\theta$ then extends to a surjective homomorphism $\theta: \btplus \ra \Cp$, which
sends $\sum p^k [x_k]$ to $\sum p^k x_k^{(0)}$. Let $[\eps]$ be the Teichm{\"u}ller lift of
$\eps$, $[\eps_1]=[(\eps^{(1)},\cdots)]$ so that $\eps_1^p=\eps$, and let
$\omega=([\eps]-1)/([\eps_1]-1)$. Then
$\theta(\omega)=1+\eps^{(1)}+\cdots+(\eps^{(1)})^{p-1}=0$, and actually,
the kernel of $\theta$ is the ideal generated by $\omega$.

Here is a simple proof: obviously, the kernel of $\theta: \etplus \ra \OO_{\Cp}/p$
is the ideal of $x \in \etplus$ such that $\vale(x) \geq 1$. Let $y$ be any element of
$\atplus$ killed by $\theta$ whose reduction modulo $p$ satisfies $\vale(\overline{y})=1$.
The map $y \atplus \ra \ker(\theta)$ is then injective, and surjective modulo $p$; since
both sides are complete for the $p$-adic topology, it is an isomorphism. Now, we just need
to observe that the element $\omega$ is killed by $\theta$ 
and that $\vale(\overline{\omega})=1$.

\begin{bibli}
These constructions are given in Fontaine's \cite{Bu88per}, 
but the reader should be warned
that the notation is rather different; for example, 
$\etplus$ is Fontaine's $\mathcal{R}$ and $\atplus$ is his
$\mathbf{A}_{\mathrm{inf}}$. 
In \cite{Bu88per}, the title of this paragraph is also explained (the
pair $(\bplus,\theta:\bplus \ra \Cp)$ is the solution of a universal problem).
The most up-to-date place to read about these rings 
is Colmez's \cite[\S 8]{Co00}.
\end{bibli}

\subsubsection{Construction of $\bdR$}\label{bdr}
Using this we can finally define $\bdR$; let $\bdR^+$ be the ring obtained by completing
$\btplus$ for the $\ker(\theta)$-adic topology, so that
$\bdR^+=\projlim_n \btplus / (\ker(\theta))^n$. In particular, every element $x \in \bdR^+$
can be written (in many ways) as a sum $x=\sum_{n=0}^{+\infty} x_n \omega^n$ with $x_n \in
\btplus$. The ring $\bdR^+$ is then naturally a $F$-algebra. Let us construct an interesting
element of this ring; since $\theta(1-[\eps])=0$, the element
$1-[\eps]$ is ``small'' 
for the topology of $\bdR^+$ and the
following series \[ -\sum_{n = 1}^{+\infty} \frac{(1-[\eps])^n}{n} \] will converge in
$\bdR^+$, to an element which we call $t$. Of course, one should think of $t$ as
$t=\log([\eps])$. For instance, if $g \in G_F$, then 
\[ g(t)=g(\log([\eps]))=\log([g(\eps^{(0)},\eps^{(1)},\cdots)])=
\log([\eps^{\chi(g)}])=\chi(g)t \]
so that $t$ is a period of the cyclotomic character. 

We now set $\bdR=\bdR^+[1/t]$, which is a field that 
we endow with the filtration defined by $\on{Fil}^i \bdR = t^i
\bdR^+$. This is the natural filtration on $\bdR$ coming from the fact
that $\bdR^+$ is a complete discrete valuation ring.
By functoriality, all the rings we have defined are equipped with a continuous
linear action of $G_K$. One can show that $\bdR^{G_K}=K$, so that if $V$ is a $p$-adic
representation, then $\ddR(V)=(\bdR \otimes_{\Qp} V)^{G_K}$ is naturally a 
filtered $K$-vector space. We say that $V$ is 
\emph{de Rham} if $\dim_K \ddR(V) = d$.

We see that $\on{Gr} \bdR \simeq \oplus_{i \in \ZZ} \Cp(i)$, and therefore, if $V$ is a de
Rham representation (a $\bdR$-admissible representation), then 
there exist $d$ integers $h_1,\cdots,h_d$ such that $\Cp \otimes_{\Qp} V
\simeq \oplus_{j=1}^d \Cp(h_j)$. A de Rham representation is therefore Hodge-Tate.
Furthermore, one sees easily that the jumps of 
the filtration on $\ddR(V)$
are precisely the opposites of
Hodge-Tate weights of $V$ (that is, $\on{Fil}^{-h_j}(D) \neq \on{Fil}^{-h_j+1}(D)$).
For example, if $V=\Qp \otimes_{\Zp} T_p E$, 
then the Hodge-Tate weights of $V$ are $0$ and $1$.

\begin{bibli}
References for this paragraph are Fontaine's 
\cite{Bu88per} for the original construction of
$\bdR$, and Colmez's \cite[\S 8]{Co00} for a more general presentation. For the behavior
of $\bdR$ under the action of some closed subgroups of $G_K$, 
one can see Iovita-Zaharescu's \cite{IZ98,IZ99}.
\end{bibli}

\subsubsection{Sen's theory for $\bdR^+$}\label{bdrsen}
Fontaine has done the analogue of Sen's theory for $\bdR^+$, that is, he defined a
$K_{\infty}[[t]]$-module $\ddif^+(V)$ which is the union of the finite dimensional 
$K_{\infty}[[t]]$-submodules of $(\bdR^+ \otimes_{\Qp} V)^{H_K}$ which are stable by
$\Gamma_K$. He then proved that $\ddif^+(V)$ is a $d$-dimensional 
$K_{\infty}[[t]]$-module endowed with a residual action of $\Gamma_K$. 
The action of $\on{Lie}(\Gamma_K)$ gives rise to a
differential operator $\nabla_V$. 
The representation $V$ is de Rham if and only if $\nabla_V$ is trivial
on $K_{\infty}((t)) \otimes_{K_{\infty}[[t]]} \ddif^+(V)$. Furthermore, one recovers
$(\dsen(V),\Theta_V)$ from $(\ddif^+(V),\nabla_V)$ simply by applying the 
map $\theta : \bdR^+ \ra \Cp$.

\begin{bibli}
This construction is carried out in Fontaine's course \cite[\S 3,4]{Fo00}, where 
$\bdR$-representations are classified.
\end{bibli}

\Subsection{The rings $\bcris$ and $\bst$}
\subsubsection{Construction of $\bcris$}\label{bcris}
One unfortunate feature of $\bdR^+$ is that it is too coarse a ring: there is
no natural extension of the natural Frobenius
$\phi:\btplus \ra \btplus$ to a continuous map $\phi: \bdR^+ \ra
\bdR^+$. For example, $\theta([\tilde{p}^{1/p}]-p) \neq 0$, so that 
$[\tilde{p}^{1/p}]-p$ is invertible in $\bdR^+$, and so $1/([\tilde{p}^{1/p}]-p) \in
\bdR^+$. But if
$\phi$ is a natural extension of  $\phi:\btplus \ra \btplus$, then one should have
$\phi(1/([\tilde{p}^{1/p}]-p))=1/([\tilde{p}]-p)$, and since $\theta([\tilde{p}]-p)=0$, 
$1/([\tilde{p}]-p) \notin \bdR^+$. 

Another way to see this is that since $\bdR^{G_L}=L$ for every finite extension $L/K$, the
existence of a canonical
Frobenius map $\phi: \bdR \ra \bdR$ would imply the existence of a Frobenius
map $\phi : \overline{K} \ra \overline{K}$, which is of course not the case.
One would still like to have a Frobenius map, and there is a natural
way to complete $\btplus$ (where one avoids inverting
elements like $[\tilde{p}^{1/p}]-p$) such that the completion is still
endowed with a Frobenius map.

The ring $\bcris^+$ is a subring of $\bdR^+$, consisting of the limits
of sequences  of elements of $\bdR^+$ which satisfy some growth condition. For example,
$\sum_{n=0}^{+\infty} p^{-n^2}t^n$ converges in $\bdR^+$ but not in $\bcris^+$.
The ring $\bcris^+$ is then equipped with a continuous Frobenius $\phi$. More precisely, 
recall that every element $x \in \bdR^+$ can be written (in many ways) as
$x=\sum_{n=0}^{+\infty} x_n \omega^n$ with $x_n \in \btplus$. One then has:
\[ \bcris^+ = \{ x \in \bdR^+,\ x=\sum_{n=0}^{+\infty} x_n \frac{\omega^n}{n!},\ 
\text{where $x_n \ra 0$ in $\btplus$} \} \]

Let $\bcris=\bcris^+[1/t]$ (note that $\bcris$ is not a field. For example,
$\omega-p$ is in $\bcris$ but not its inverse - this is not so easy to see); one can show
that $(\bcris)^{G_K}=F$. Those representations $V$ of $G_K$ 
which are $\bcris$-admissible are called
\emph{crystalline}, and using Fontaine's construction 
one can therefore associate to every such $V$
a filtered
$\phi$-module
$\dcris(V)=(\bcris \otimes_{\Qp} V)^{G_K}$ 
(a \emph{filtered $\phi$-module} $D$ is an $F$-vector space with a decreasing,
exhaustive and separated filtration indexed by $\ZZ$
on $K \otimes_F D$, 
and a $\sigma_F$ semi-linear map $\phi: D \ra D$. We
do not impose any compatibility condition between $\phi$ and $\on{Fil}$). One can associate
to a filtered
$\phi$-module $D$ two polygons: its Hodge polygon $P_H(D)$, coming from the filtration, and
its Newton polygon $P_N(D)$, coming from the slopes of $\phi$. We say that $D$ is
\emph{admissible} 
if for every subobject $D'$ of $D$,
the Hodge polygon of $D'$ lies below the Newton polygon of $D'$, and
the endpoints 
of the Hodge and Newton polygons of $D$
are the same. 
One can show that $\dcris(V)$ is always admissible.

Furthermore, a theorem of Colmez and Fontaine shows that the functor $V \mapsto \dcris(V)$ is
an equivalence of categories between the category of crystalline representations 
and the category of admissible filtered $\phi$-modules
\footnote{admissible modules were previously called \emph{weakly admissible}, but
  since Colmez and Fontaine showed that being weakly admissible is the
  same as being admissible (previously, $D$ was said to be admissible if there
  exists some $V$ such that $D=\dcris(V)$), we can drop the ``weakly''
  altogether.}. 

\begin{bibli}
The construction of $\bcris$ can be found in 
Fontaine's \cite{Bu88per} or Colmez's \cite[\S 8]{Co00}. One should also
look at Fontaine's \cite{Bu88sst} for information on filtered $\phi$-modules. 
The theorem of Colmez-Fontaine is proved in Colmez-Fontaine's \cite{CoFo}, as well as in
Colmez's \cite[\S 10]{Co00} and it is reviewed in
Fontaine's \cite[\S 5]{Fo00}. The ring $\bcris^+$ has an
interpretation in crystalline cohomology, see Fontaine's \cite{Fo82}
and Fontaine-Messing's \cite{FMsg}.
\end{bibli}

\subsubsection{Example: elliptic curves}\label{polyg}
If $V=\Qp \otimes_{\Zp} T_pE$, where $E$ is an
elliptic curve over $F$
with good ordinary reduction, then $\dcris(V)$ is a $2$-dimensional
$F$-vector space with a basis
$x,y$, and there exists
$\lambda \in F$ and $\alpha_0, \beta_0 \in \OO_F^*$ depending on $E$ such that:
\[  \begin{cases}
\phi(x) & = \alpha_0 p^{-1}x \\
\phi(y) & = \beta_0 y
\end{cases}
\qquad\text{and}\qquad  \on{Fil}^i \dcris(V) =
\begin{cases}
\dcris(V) & \text{if $i \leq -1$} \\
(y+\lambda x)F  & \text{if $i = 0$} \\
\{0\}  & \text{if $i \geq 1$}
\end{cases} \]

The Newton and Hodge polygons of $\dcris(V)$ are then as follows:
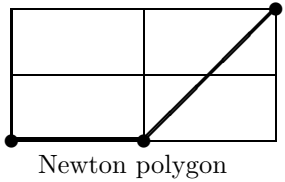
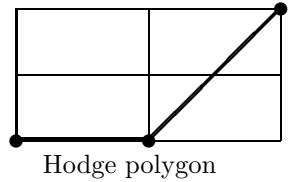
\begin{center}
\begin{picture}(300,75)(0,0)
\thicklines
\put(0,20){\line(1,0){50}}
\put(50,20){\line(1,1){50}}

\put(200,20){\line(1,0){50}}
\put(250,20){\line(1,1){50}}

\put(0,21){\line(1,0){50}}
\put(50,21){\line(1,1){50}}

\put(200,21){\line(1,0){50}}
\put(250,21){\line(1,1){50}}

\thinlines
\put(0,20){\line(0,1){50}}
\put(0,20){\line(1,0){100}}

\put(0,45){\line(1,0){100}}
\put(50,20){\line(0,1){50}}

\put(100,70){\line(0,-1){50}}
\put(100,70){\line(-1,0){100}}

\put(200,20){\line(0,1){50}}
\put(200,20){\line(1,0){100}}

\put(200,45){\line(1,0){100}}
\put(250,20){\line(0,1){50}}

\put(300,70){\line(0,-1){50}}
\put(300,70){\line(-1,0){100}}

\put(0,20){\circle*{5}}
\put(50,20){\circle*{5}}
\put(100,70){\circle*{5}}
\put(200,20){\circle*{5}}
\put(250,20){\circle*{5}}
\put(300,70){\circle*{5}}

\put(10,8){Newton polygon}
\put(210,8){Hodge polygon}
\end{picture}
\end{center}

If on the other hand, an elliptic curve
$E$ has good supersingular reduction, then the operator
$\phi:\dcris(V) \ra \dcris(V)$ is irreducible and the Newton 
and Hodge polygons are as follows:
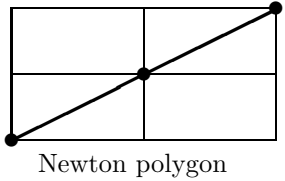
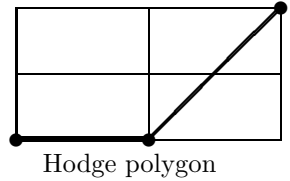
\begin{center}
\begin{picture}(300,75)(0,0)
\thicklines
\put(0,20){\line(2,1){100}}

\put(200,20){\line(1,0){50}}
\put(250,20){\line(1,1){50}}

\put(1,20){\line(2,1){100}}

\put(200,21){\line(1,0){50}}
\put(250,21){\line(1,1){50}}

\thinlines
\put(0,20){\line(0,1){50}}
\put(0,20){\line(1,0){100}}

\put(0,45){\line(1,0){100}}
\put(50,20){\line(0,1){50}}

\put(100,70){\line(0,-1){50}}
\put(100,70){\line(-1,0){100}}

\put(200,20){\line(0,1){50}}
\put(200,20){\line(1,0){100}}

\put(200,45){\line(1,0){100}}
\put(250,20){\line(0,1){50}}

\put(300,70){\line(0,-1){50}}
\put(300,70){\line(-1,0){100}}

\put(0,20){\circle*{5}}
\put(50,45){\circle*{5}}
\put(100,70){\circle*{5}}
\put(200,20){\circle*{5}}
\put(250,20){\circle*{5}}
\put(300,70){\circle*{5}}

\put(10,8){Newton polygon}
\put(210,8){Hodge polygon}
\end{picture}
\end{center}

In both cases, it is clear that $\dcris(V)$ is admissible.

\begin{bibli}
For basic facts about elliptic curves, see for example Silverman's
\cite{Si86,Si96}. For basic facts on Newton polygons, see the first
chapter of \cite{DGS}.
\end{bibli}

\subsubsection{Semi-stable representations}\label{bst}
If an elliptic curve
$E$ has bad semi-stable reduction, then $V$ is not 
crystalline but it is \emph{semi-stable},
that is, it is $\bst$-admissible where $\bst=\bcris[Y]$, where we have decided that 
$\phi(Y)=Y^p$ and $g(Y)=Y+c(g)t$, where $c(g)$ is defined by the formula
$g(p^{1/p^n})=p^{1/p^n} (\eps^{(n)})^{c(g)}$. Of course, the definition of $Y$ depends on a
number of choices, but two such
$\bst$'s are isomorphic. In addition to a Frobenius, $\bst$ is
equipped with the \emph{monodromy map} $N=-d/dY$.

Let $\tilde{p} \in \etplus$ be an element such that
$\tilde{p}^{(0)}=p$, and 
let $\log[\tilde{p}] \in \bdR^+$ be the element defined by
\[ \log[\tilde{p}] = \log_p(p) - \sum_{n=1}^{+\infty} \frac{(1-[\tilde{p}]/p)^{n-1}}{n}.  \]
One can define a Galois equivariant and $\bcris$-linear embedding of $\bst$ into $\bdR$, by
mapping $Y$ to 
$\log[\tilde{p}] \in \bdR^+$, 
but doing so requires us to make a choice of $\log_p(p)$. As a
consequence, there is no canonically defined filtration on
$\dst(V)$, only on $\ddR(V)$: one has to be a little careful about this. This in contrast to
the fact that the inclusion $K \otimes_F \dcris(V) \subset \ddR(V)$ is canonical. It is
customary to choose $\log_p(p)=0$ which is what we'll always assume from
now on.

One can then associate to every semi-stable representation $V$ a filtered $(\phi,N)$-module
and Colmez and Fontaine showed that the functor $V \mapsto \dst(V)$ is an equivalence of
categories between the category of semi-stable representations and the category of
admissible filtered $(\phi,N)$-modules.

\begin{bibli}
See the references for paragraph \ref{bcris} on $\bcris$. See
\ref{pereq} for Fontaine's original definition of $\bst$.
\end{bibli}

\subsubsection{Frobenius and filtration}\label{exseq}
Although the ring $\bcris$ is endowed with both a Frobenius map and the filtration
induced by $\bcris \subset \bdR$, these two structures have little compatibility. For
example, here is a funny exercise: let $r=\{r_n\}_{n \geq 0}$ be a sequence with $r_n \in
\ZZ$. Show that there exists an element $x_r \in \bcris$ (in the image of $\phi^n$ for all
$n \geq 0$)
such that for every $n \geq 0$, one
has $\phi^{-n}(x_r) \in \on{Fil}^{r_n} \bdR \setminus \on{Fil}^{r_n+1} \bdR$ (for a
solution, see paragraph \ref{algfun}).  The reader
should also be warned that $\bcris^+ \subset \on{Fil}^0 \bcris
= \bcris \cap \bdR^+$ but that the latter space is
much larger. It is true however that if $\bcris'$ is the set of elements $x \in \bcris$ such
that for every $n \geq 0$, one has $\phi^n(x) \in \on{Fil}^0 \bcris$,
then $\phi(\bcris') \subset \bcris^+ \subset \bcris'$
($\phi^2(\bcris')$ if $p=2$).

Given the above facts, it is rather surprising that there is a relation of some sort between
$\phi$ and $\on{Fil}$. One can show that the natural map $\bcris^{\phi=1} \ra
\bdR/\bdR^+$ is surjective, and that its kernel is $\Qp$. 
This gives rise to an exact sequence
\[ 0 \ra \Qp \ra \bcris^{\phi=1} \ra \bdR/\bdR^+ \ra 0 \]
called the \emph{fundamental exact sequence}. 
It will be used later on to define Bloch-Kato's
exponential.

\begin{bibli}
See Fontaine's \cite{Bu88per} and \cite{Bu88sst} or Colmez's \cite[\S 8]{Co00}. For
Bloch-Kato's exponential, see paragraph \ref{expbk} as well as 
Bloch-Kato's \cite[\S 3]{BK91} and Kato's \cite{Ka93}.
\end{bibli}

\subsubsection{Some remarks on topology}\label{topo}
We'll end this section with a few remarks 
on the topologies of the rings we just introduced.
Although $\bdR^+$ is a discretely valued ring, and complete for that
valuation, the natural
topology on $\bdR^+$
is weaker. It is actually the topology of the projective limit on
$\bdR^+ = \projlim_n \btplus/\ker(\theta)^n$, and the topology of $\btplus=\atplus[1/p]$
combines the $p$-adic topology and the topology of the residue ring $\atplus/p=\etplus$
which is a valued ring. In particular, $\bdR^+/\ker(\theta)^n$ is
$p$-adic Banach space, which makes $\bdR^+$ into a $p$-adic Fr\'echet space.

The topology on $\bcris$ is quite unpleasant, as Colmez points out: ``By the construction of
$\bcris^+$, the sequence $x_n=\omega^{p^n-1}/(p^n-1)!$ does not converge to $0$ in
$\bcris^+$, but the sequence $\omega x_n$ does; we deduce from this the fact that the
sequence $t x_n$ converges to $0$ in $\bcris^+$, and therefore
that $x_n \ra 0$ in $\bcris$, so
that the topology of $\bcris^+$ induced by that of $\bcris$ is not the natural topology of
$\bcris^+$.'' The reason for this is that the sequence $n!$ converges to $0$ in a pretty
chaotic manner, and it is more convenient to use the ring 
\[ \bmax^+ = \{ x \in \bdR^+,\ x=\sum_{n=0}^{+\infty} x_n \frac{\omega^n}{p^n},\ 
\text{where $x_n \ra 0$ in $\btplus$} \}, \]
which is also endowed with a Frobenius map.
In any case, the periods of crystalline representations live in
\[ \btrigplus{}[1/t]=\cap_{n=0}^{+\infty} \phi^n \bcris^+ [1/t]= \cap_{n=0}^{+\infty} \phi^n
\bmax^+[1/t] \] 
because they live in finite dimensional $F$-vector subspaces of $\bcris$ stable by $\phi$.

Finally, let us mention an interesting result of Colmez, that has yet
to be applied: 
$\Kbar$ is naturally a subring of $\bdR^+$, and he
showed that
$\bdR^+$ is the completion of
$\Kbar$ for the induced topology, which is finer than the $p$-adic
topology 
(meaning there are
more open sets). 
This generalizes an earlier result of Fontaine, who showed that
$\Kbar$ is dense $\bdR^+/t^2$.
The topology of $\Kbar$ induced by $\bdR^+$ is a bit like the
``uniform convergence of a function and all its derivatives'', if one
views $x \in \Kbar$ as an algebraic function of $p$. 
For example, the series $\sum_{n=0}^{+\infty} p^n \eps^{(n)}$ is not
convergent in
$\bdR^+$. A series which converges in $\bdR^+$ does so in $\Cp$, so we get a map
$\theta: \bdR^+ \ra \Cp$, which coincides with the one previously defined.

\begin{bibli}
The remark on the topology of $\bcris$ can be found in Colmez's \cite[III]{Co98}, and
Colmez's theorem is proved in the appendix to Fontaine's
\cite{Bu88per}. 
Fontaine's earlier result was used by Fontaine and Messing in
\cite{FMsg}. The ring
$\btrigplus{}$ has an interpretation in rigid cohomology, as was explained to me by
Berthelot in \cite{Bt01}.
\end{bibli}

\Subsection{Application: Tate's elliptic curve}
We will now explicitly show that if $E$ is an elliptic 
curve with bad semi-stable reduction, then $V=\Qp \otimes_{\Zp} T_p E$ is $\bdR$-admissible.
After that, we will show that $V$ is actually semi-stable. We'll assume throughout this
section that $K=F$. 

\subsubsection{Tate's elliptic curve}\label{tatec}
Let $q$ be a formal parameter and define
\[ s_k(q) = \sum_{n=1}^{+\infty} \frac{n^k q^n}{1-q^n} \qquad
a_4(q)=-s_3(q) \qquad
a_6(q) = - \frac{5 s_3(q)+7 s_5(q)}{12} \]
\[ x(q,v)=\sum_{n=-\infty}^{+\infty} \frac{q^n v}{(1-q^n v)^2}-2 s_1(q) \qquad
y(q,v)=\sum_{n=-\infty}^{+\infty} \frac{(q^n v)^2}{(1-q^n v)^3}+ s_1(q).
\]
All those series are convergent if $q \in p \OO_F$ and 
$v \notin q^{\ZZ} = \langle q \rangle$ (the multiplicative subgroup of
$F^*$ generated by $q$). 
For such $q\neq 0$, let $E_q$ be the elliptic
curve defined by the equation $y^2+xy=x^3 + a_4(q) x + a_6(q)$. The theorem of Tate
is then: the elliptic curve $E_q$ is defined over $F$, it has bad semi-stable
reduction, and $E_q$ is uniformized by $\Fbar^*$, that is, there
exists a map
$\alpha:\Fbar^* \ra E_q(\Fbar)$, given by
\[ v \mapsto \begin{cases} 
(x(q,v),y(q,v)) & \text{ if } v \notin q^{\ZZ} \\
0  & \text{ if } v \in q^{\ZZ} \end{cases} \]
which induces an isomorphism of groups with $G_F$-action $\Fbar^*/\langle q \rangle \ra
E_q(\Fbar)$.

Furthermore, if $E$ is an elliptic curve over $F$ with bad semi-stable reduction, then
there exists $q$ such that $E$ is isomorphic to $E_q$ over $F$.

\begin{bibli}
For basic facts about Tate's elliptic curve, see Silverman's \cite[V.3]{Si96}.
\end{bibli}

\subsubsection{The $p$-adic representation attached to $E_q$}\label{ordi}
Using Tate's theorem, we can give an explicit description of $T_p(E_q)$. Let $\eps^{(i)}$ be
the $p^i$-th roots of unity chosen in \ref{nota}
and let $q^{(i)}$ be elements defined by
$q^{(0)}=q$ and the requirement that $(q^{(i+1)})^p=q^{(i)}$. Then $\alpha$ induces
isomorphisms
\[ \begin{CD}
\Fbar^* /\langle q \rangle @>>> E_q(\Fbar) \\
\{ x \in \Fbar^* /\langle q \rangle, x^{p^n} \in \langle q \rangle \} @>>>
E_q(\Fbar)[p^n] 
\end{CD} \]
and one sees that $\{ x \in \Fbar^* /\langle q \rangle, x^{p^n} \in \langle q \rangle \}
=\{ (\eps^{(n)})^i (q^{(n)})^j, 0 \leq i,j < p^n-1 \}$. The elements $\eps^{(n)}$ and
$q^{(n)}$ therefore form a basis of $E_q(\Fbar)[p^n]$, so that a basis of $T_p(E_q)$ is given
by $e=\projlim_n \eps^{(n)}$ and $f=\projlim_n q^{(n)}$. This makes it possible to compute
explicitly the Galois action on $T_p(E_q)$. We have $g(e)=\projlim_n g(\eps^{(n)})=\chi(g)e$
and $g(f)=\projlim_n g(q^{(n)}) = \projlim_n q^{(n)} (\eps^{(n)})^{c(g)}=f+c(g)e$ where $c(g)$
is some $p$-adic integer, determined by the fact that $g(q^{(n)})=q^{(n)}(\eps^{(n)})^{c(g)}$.
Note that $[g \mapsto c(g)] \in H^1(F,\Zp(1))$.
The matrix of $g$ in the basis $(e,f)$ is therefore given by 
\[ \begin{pmatrix}
\chi(g) & c(g) \\
0 & 1 \end{pmatrix} \]

\subsubsection{$p$-adic periods of $E_q$}\label{pereq}
We are looking for $p$-adic periods of $V=\Qp \otimes_{\Zp} T_p(E_q)$ which live in $\bdR$,
that is for elements of $(\bdR \otimes_{\Qp} V)^{G_F}$. An obvious candidate is $t^{-1}
\otimes e$ since
$g(t)=\chi(g) t$ and $g(e) = \chi(g) e$. Let us look for a second element of $(\bdR
\otimes_{\Qp} V)^{G_F}$, of the form $a \otimes e + 1 \otimes f$. We see that this element
will be fixed by $G_F$ if and only if $g(a)\chi(g)+c(g)=a$. 

Let $\tilde{q}$ be the element of $\etplus$ defined by $\tilde{q}=(q^{(0)},q^{(1)},\cdots)$.
Observe that we have $g(\tilde{q}) = (g(q^{(0)}),g(q^{(1)}),\cdots) = \tilde{q} \eps^{c(g)}$,
and that $\theta([\tilde{q}]/q^{(0)}-1)=0$, so that $[\tilde{q}]/q^{(0)}-1$ is small in the
$\ker(\theta)$-adic topology. The series
\[ \log_p(q^{(0)}) - \sum_{n=1}^{+\infty}
\frac{(1-[\tilde{q}]/q^{(0)})^n}{n} \]
therefore converges in $\bdR^+$ to an element which we call $u$. One should think of $u$ as
being $u=\log([\tilde{q}])$. In particular, $g(u)=g(\log([\tilde{q}]))=
\log([g(\tilde{q})])=\log([\tilde{q}])+c(g)\log([\eps])=u+c(g)t$, and we readily see that
$a=-u/t$ satisfies the equation $g(a)\chi(g)+c(g)=a$. A basis of $\ddR(V)=
(\bdR \otimes_{\Qp} V)^{G_F}$ is therefore given by 
\[ \begin{cases}
x = t^{-1} \otimes e \\
y = -ut^{-1} \otimes e + 1 \otimes f 
\end{cases} \]
and this shows that $T_p(E_q)$ is $\bdR$-admissible. 
Furthermore, one sees that $\theta(u-\log_p(q^{(0)}))=0$, so that
$u-\log_p(q^{(0)})$ is divisible by $t$ and 
\[ \on{Fil}^i \ddR(V) = \begin{cases}
 \ddR(V) & \text{ if } i \leq -1 \\
(y+\log_p(q^{(0)})x)F & \text{ if } i = 0 \\
 \{ 0 \} & \text{ if } i \geq 1 \\
\end{cases} \]

This gives us a description of $\ddR(V)$. We shall now prove that $V$ is semi-stable. It's
clearly enough to show that $t, u \in \bst^+$. The series which defines $t$ converges in
$\bcris^+$ (that is, the cyclotomic character is crystalline), and the series which defines
$\log[\tilde{q}/\tilde{p}^{v_p(q)}]$ also does. As a consequence, one can write 
$u=v_p(q)Y+\log[\tilde{q}/\tilde{p}^{v_p(q)}] \in \bst^+$. This implies that $V$ is
semi-stable. Actually, Fontaine defined $\bst$ so that it would
contain $\bcris$ and a period of $E_q$, so that the computation of
this paragraph is a little circular.

Let us compute the action of Frobenius in the case of Tate's elliptic curve. On a ring of
characteristic $p$, one expects Frobenius to be $x \mapsto x^p$, and therefore $\phi([x])$
should be $[x^p]$ so that $\phi(\log[x])=p\log[x]$. In particular, one has $\phi(t)=pt$
and $\phi(u)=pu$ and the action of Frobenius on $\dst(V)$ is therefore given by 
$\phi(x)=p^{-1}x$ and $\phi(y)= y$. Let us point out one more time that the filtration is
defined on $\ddR(V)$, and that the identification $\dst(V) \simeq \ddR(V)$ depends on
a choice of $\log_p(p)$. The $p$-adic number
$\log_p(q^{(0)}/p^{v_p(q^{(0)})})$ is canonically attached to $V$ and
is called the \emph{$\ell$-invariant} of $V$.

\subsubsection{Remark: Kummer theory}\label{kummer}
What we have done for Tate's elliptic curve is really a consequence of the fact that 
$V=\Qp \otimes_{\Zp} T_p E_q$ is an extension of $\Qp$ by $\Qp(1)$, namely that there is an
exact sequence $0 \ra \Qp(1) \ra V \ra \Qp \ra 0$. All of these extensions are classified by
the cohomology group $H^1(K,\Qp(1))$, which is described by Kummer theory. Recall that
for every $n \geq 1$, there is an 
isomorphism $\delta_n : K^*/(K^*)^{p^n} \ra H^1(K,\mu_{p^n})$. 
By taking the projective limit over $n$, we get a map $\delta: \hat{K^*} \ra 
H^1(K,\Zp(1))$ because $\projlim_n \mu_{p^n} \simeq \Zp(1)$ once we have chosen a
compatible sequence of $\eps^{(n)}$.  
By tensoring with $\Qp$, we get an isomorphism
$\delta: \Qp \otimes_{\Zp} \hat{K^*} \ra  H^1(K,\Qp(1))$ 
which is defined in the following
way: if $q = q^{(0)} \in \Qp \otimes_{\Zp} \hat{K^*}$, 
choose a sequence $q^{(n)}$ such that 
$(q^{(n)})^p=q^{(n-1)}$, and define 
$c=\delta(q)$ by $(\eps^{(n)})^{c(g)}=g(q^{(n)})/(q^{(n)})$.
Of course, this depends on the choice of $q^{(n)}$, but two different choices give
cohomologous cocycles.

It is now easy to show that every extension of $\Qp$ by $\Qp(1)$ is semi-stable. This is
because $t c(g)=g(\log[\tilde{q}])-\log[\tilde{q}]$ 
with notations similar to those above, and $\tilde{q}=(q^{(n)})$. If 
$q\in \Qp \otimes_{\Zp} \hat{\OO_K^*}$ then the series which defines $\log[\tilde{q}]$
converges in $\bcris^+$ and the extension $V$ is crystalline. In
general, 
if $q\in \Qp \otimes_{\Zp} \hat{K^*}$, then
$\log[\tilde{q}]$
will be in $\bcris^++v_p(q)Y \subset \bst^+$. The $F$-vector space $\dst(V)$ will then have
a basis $x=t^{-1} \otimes e$ and $y=-\log[\tilde{q}]t^{-1} \otimes e + 1 \otimes f$ so that
$\phi(x)=p^{-1}x$ and $\phi(y)=y$. If one chooses $\log_p(p)=0$, then the filtration on
$\ddR(V)$ is given by $\on{Fil}^0 \ddR(V)  = (y + \log_p(q)x)F$.

\Subsection{$p$-adic representations and Arithmetic Geometry}
\subsubsection{Comparison theorems}\label{tsuji}
If $X/K$ is a proper smooth variety over $K$, then
by a \emph{comparison theorem}, we mean a theorem relating 
$H^i_{\mathrm{\acute{e}t}}(X_{\overline{K}},\Qp)$ 
and $H^i_{\mathrm{dR}}(X/K)$.

It was shown early on by Fontaine that the Tate modules $V=\Qp
\otimes_{\Zp} T_pA$ of all abelian varieties $A$ are de Rham
(he actually showed in a letter to Jannsen that
they were potentially semi-stable), 
and that $\ddR(V)$ is
isomorphic to the dual of the de Rham cohomology of $A$. 
Fontaine and Messing then found another proof, in which they
explicitly construct a pairing between $V$ (interpreted as a quotient of the \'etale
$\pi_1(A)$) and $H^1_{\mathrm{dR}}(A/K)$ 
(interpreted as the group of isomorphism classes of
vectorial extensions of $A$). One should remember that for an abelian variety $A$, we have
$\on{Hom}_{G_K}(T_pA,\Zp(1)) \simeq H^1_{\mathrm{\acute{e}t}}(A_{\overline{K}},\Zp)$.

After that, Fontaine and Messing proved the comparison theorem for the 
$H^i_{\mathrm{\acute{e}t}}(X_{\overline{K}},\Qp)$ of proper smooth $X$
for $i \leq p-1$ and $K/\Qp$ finite unramified. These results were
then extended by Kato and his school (Hyodo, Tsuji). 
Finally, the general statement that for a variety $X/K$, one can recover 
$H^i_{\mathrm{dR}}(X/K)$ from the 
data of $V=H^i_{\mathrm{\acute{e}t}}(X_{\overline{K}},\Qp)$ as
a $p$-adic representation was shown by Tsuji. He  
showed that if $X$ has semi-stable reduction, then 
$V=H^i_{\mathrm{\acute{e}t}}(X_{\overline{K}},\Qp)$ is $\bst$-admissible. 
A different proof was given by Niziol (in the good reduction case) and
also by Faltings (who proved that $V$ is crystalline if $X$ has good
reduction and that $V$ is de Rham otherwise).

In the case of an abelian variety, 
the rings $\bcris$ and $\bst$ are exactly what it takes to decide, from the
data of $V$ alone, whether $A$ has good or semi-stable reduction. Indeed, Coleman-Iovita and
Breuil showed that $A$ has good reduction if and only if $V$ is crystalline, and that $A$
has semi-stable reduction if and only if $V$ is semi-stable. This 
can be seen as a $p$-adic analogue
of the ($\ell$-adic) N\'eron-Ogg-Shafarevich criterion. 

In another direction, Fontaine and Mazur have conjectured the following: 
let $V$ be a $p$-adic
representation of $\on{Gal}(\overline{\QQ}/L)$ where $L$ is a finite
extension of $\QQ$. Then, $V$ ``comes from geometry'' if and only if it is 
unramified at all
but finitely many primes $\ell$, and if its restriction to all decomposition groups at $p$
is  potentially semi-stable.

Note that among all potentially semi-stable representations $V$ of $G_K$, where $K$ is a
$p$-adic field, there are many which do not come from geometry: indeed, 
if $V=H^i_{\mathrm{\acute{e}t}}(X_{\overline{K}},\Qp)$ then the eigenvalues of
$\phi$ on $\dst(V)$ should at least be Weil numbers.

\begin{bibli}
There were many partial results before Tsuji's theorem was proved in
\cite{Ts99} (see \cite{Ts02} for a survey), 
and we refer the reader to the bibliography of that article. 
For a different approach (integrating forms on cycles), 
see Colmez's \cite{CoAs}. 

The conjecture of Fontaine and Mazur
was proposed by them in \cite{FM93}. 
There is little known in that direction, but there is some
evidence in dimension $2$: see Taylor's \cite{RT01}. Kisin has some
interesting results too in \cite{Ki02}.

Regarding the criteria for good or semi-stable reduction, see Coleman-Iovita's
\cite{CI99} and Breuil's \cite{Br99}.
\end{bibli}

\subsubsection{Weil-Deligne representations}\label{weild}
Let $V$ be a potentially semi-stable representation of $G_K$, so that
there exists $L$, a finite
extension of $K$ such that the restriction of $V$ to $G_L$ is
semi-stable. One can then consider the $F$-vector space $\dst^L(V)=(\bst \otimes_{\Qp}
V)^{G_L}$. It is a finite dimensional $(\phi,N)$-module with an action of $\on{Gal}(L/K)$.
One can attach to such an object several interesting invariants: $L$-factors,
$\epsilon$-factors, and a representation of the Weil-Deligne group. 

In particular, if $E$ is an elliptic curve, 
one can recover from the $p$-adic representation $T_p E$ pretty much the
same information as from the $\ell$-adic representation $T_{\ell} E$.

\begin{bibli}
The action of the Weil-Deligne group on $\dst^L(V)$ was defined by
Fontaine in \cite{Bu88ladic}.
\end{bibli}

\[ \clubsuit\diamondsuit\spadesuit\heartsuit \]

\newpage

\section{Fontaine's $(\phi,\Gamma)$-modules}
\Subsection{The characteristic $p$ theory}
A powerful tool for studying $p$-adic representations is Fontaine's theory of
$(\phi,\Gamma)$-modules. We will first define $\phi$-modules for representations of 
the Galois group of a local field of characteristic $p$ (namely $k((\pi))$) and then apply
this to the characteristic zero case, making 
use of Fontaine-Wintenberger's theory of the field
of norms. To simplify the exposition, we will sometimes assume in this chapter that
$K=F$ (meaning that $K$ is absolutely unramified).

\subsubsection{Local fields of characteristic $p$}\label{charp}
Let $\pi_K$ be a formal variable (for now) and let 
$\mathbf{A}_K$ be the ring \[ \mathbf{A}_K =
\{ \sum_{k=-\infty}^{\infty}  a_k \pi_K^k,\ a_k \in \OO_F,\ a_{-k} \ra 0 \}, \]
so that $\mathbf{A}_K/p=k_K((\pi_K))$.
The ring $\mathbf{A}_K$ 
(which is an example of a Cohen ring, as in \ref{witt})
is endowed with actions of $\phi$ and $\Gamma_K$, such that
$\phi(\pi_K)=\pi_K^p \mod{p}$. 
The exact formulas depend on $K$, but 
if $K=F$ then $\phi(\pi_K)=(1+\pi_K)^p-1$ and if
$\gamma\in\Gamma_K$, then $\gamma(\pi_K)=(1+\pi_K)^{\chi(\gamma)}-1$. We won't use the
action of
$\Gamma_K$ in the ``characteristic $p$ case''. Let $\e_K=k_K((\pi_K))=\mathbf{A}_K/p$, and
let
$\e$ be the  separable closure of $\e_K$. Let $G_{E_K}$ be the Galois group of $\e/\e_K$. In
this paragraph, we will look at $p$-adic representations of $G_{E_K}$, that is, finite
dimensional
$\Qp$-vector spaces $V$, endowed with a continuous linear action of $G_{E_K}$.

Let $\mathbf{B}_K$ be the fraction field of $\mathbf{A}_K$ (one only needs to invert $p$).
A \emph{$\phi$-module} $M$ is a 
finite dimensional $\mathbf{B}_K$-vector space with a
semi-linear action of $\phi$. We say that $M$ is 
\emph{\'etale} (or \emph{slope $0$} or also
\emph{unit-root}) if $M$ admits an $\mathbf{A}_K$-lattice 
$M_A$ which is stable by $\phi$ and such
that $\phi^* M_A = M_A$. This follows for example from $\phi(M_A) \subset M_A$ and $p
\nmid \det(\phi)$. The first result is that there is an equivalence of categories 
\[ \{ \text{$p$-adic representations of $G_{E_K}$}\} \qquad
\longleftrightarrow
\qquad \{ \text{\'etale $\phi$-modules} \} \]

Let us explain where this comes from. The correspondence 
$T \mapsto (\e \otimes_{\Fp} T)^{G_{E_K}}$ is (by Hilbert 90) an equivalence 
of categories between the category of $\Fp$-representations of $G_{E_K}$, and \'etale
$\e_K$-modules. Let $\mathbf{A}$ be a Cohen ring over $\e$
(we will give a more precise definition of
$\mathbf{A}$ below. Suffice it to say that $\mathbf{A}$ should be the ring of Witt vectors
over $\e$, but $\e$ is not perfect, so that there are several possible choices for
$\mathbf{A}$). 
The ring $\mathbf{A}$ is endowed with an action of 
$G_{E_K}$ and $\mathbf{A}^{G_{E_K}}=\mathbf{A}_K$.
Then by lifting things to characteristic $0$ and inverting
$p$, we get an equivalence of categories between the category of $\Qp$-representations of
$G_{E_K}$, and \'etale $\mathbf{B}_K$-modules with a Frobenius (these
constructions are mathematical folklore and were used, for example, by
Bloch and Katz).

We will now give a construction of the ring of periods $\mathbf{A}$. Let $\etplus$ be the
ring introduced in \ref{thickcp} and let $\et$ be the field of fractions of $\etplus$. Then
$\e_K$ embeds in $\et$. For example, if $K=F$, then $\e_F = k((\eps-1)) \subset \et$. Let
$\e$ be the completion of the separable closure of $\e_K$ in $\et$. 

One can show that $\et$ is the
algebraic closure of $\e_F$ so that $\et$ is the completion of the perfection of $\e$. By a
theorem of Ax, $\et$ is actually also the completion of $\e$.

Let $\at=W(\et)$ and $\bt=\at[1/p]$. It is easy to see (at least when $K=F$) that
$\mathbf{B}_K$ is a subfield of $\bt$, with $\pi_F=[\eps]-1$. 
If $K \neq F$, then one should take for $\pi_K$ an element of
$\mathbf{A}$ whose image modulo $p$ is a uniformiser of $\e_K=\e^{H_K}$.
Let $\mathbf{B}$ be the completion of the maximal
unramified extension of $\mathbf{B}_K$ in $\bt$, and $\mathbf{A}=\mathbf{B} \cap \at$. The
field $\mathbf{B}$ is endowed with an action of $G_{E_K}$, and one indeed has
$\mathbf{B}^{G_{E_K}}=\mathbf{B}_K$. It is clear that $\mathbf{B}$ is endowed with a
Frobenius map $\phi$.

\begin{bibli}
These ideas are in the folklore, see Katz in \cite[chap 4]{Ka73}. 
We gave their local version, which is in Fontaine's \cite[A1]{Fo91}.
\end{bibli}

\subsubsection{Representations of $G_{E_K}$ and differential equations}\label{matsuda}
Let us mention an application of the theory we just sketched.
Let $\delta$ be the differential operator defined by $\delta(f(\pi))
= (1+\pi)df/d\pi$ on the field $\mathbf{B}_F$. This operator extends to $\mathbf{B}$ because
it extends to the maximal unramified extension of $\mathbf{B}_F$, and then to its closure by
continuity. One can use it to associate to every $p$-adic representation of $G_{E_K}$ a
$\mathbf{B}_F$-vector space with a Frobenius $\phi$ and a differential operator $\delta$
which satisfy $\delta\circ\phi=p\phi\circ\delta$. When the action of the
inertia of $G_{E_K}$ factors
through a finite quotient on a representation $V$, then there exists a basis of $\dfont(E)$
in which
$\delta$ is overconvergent (in the sense of \ref{overcv} below). 
One can use this fact to associate to
every potentially unramified representation of $G_{E_K}$ an overconvergent differential
equation. This condition ($\phi$ and $\delta$ overconvergent) is much
stronger than merely requiring $\phi$ to be overconvergent (which
happens very often, see \ref{overcv}).

\begin{bibli}
There are many interesting parallels between the theory of finite 
Galois representations in characteristic $p$ and differential equations: see 
Crew's \cite{RC87,RC00} and Matsuda's \cite{Ma97} for a starting point.
\end{bibli}

\Subsection{The characteristic zero theory}
\subsubsection{The field of norms}\label{fofnor}
The next step of the construction is the theory of the field of norms (of Fontaine and
Wintenberger) which gives a canonical isomorphism between $G_{E_K}$ and $H_K$. Let
$\mathcal{N}_K$ be the set $\projlim_n K_n$ where the transition maps are given by
$N_{K_n/K_{n-1}}$, so that $\mathcal{N}_K$ is the set of sequences
$(x^{(0)},x^{(1)},\cdots)$ with $x^{(n)} \in K_n$ and $N_{K_n/K_{n-1}}(x^{(n)}) =
x^{(n-1)}$. If we define a ring structure on $\mathcal{N}_K$ by
\[ (xy)^{(n)}= x^{(n)}y^{(n)} \quad\text{and}\quad
(x+y)^{(n)}=\lim_{m \ra +\infty} N_{K_{n+m}/K_n}(x^{(n+m)}+y^{(n+m)}), \]
then
$\mathcal{N}_K$ is actually a field, 
called the \emph{field of norms} of $K_{\infty}/K$. It is
naturally endowed with an action of $H_K$. Furthermore, for every finite Galois extension
$L/K$, $\mathcal{N}_L/\mathcal{N}_K$ is a finite Galois extension whose Galois group is
$\on{Gal}(L_{\infty}/K_{\infty})$, and every finite Galois extension of $\mathcal{N}_K$ is
of this kind so that the absolute
Galois group of $\mathcal{N}_K$ is naturally isomorphic to $H_K$.

On the other hand, one can prove that $\mathcal{N}_K$ is a local field of characteristic $p$
isomorphic to $\e_K \simeq k((\pi_K))$. More precisely, by ramification theory, the map
$N_{K_n/K_{n-1}}$ is close to the $p$-th power map and there is therefore a well-defined ring
homomorphism from $\mathcal{N}_K$ to $\et$ given by sending $(x^{(n)}) \in \mathcal{N}_K$ 
to $(y^{(n)}) \in \et$ where
$y^{(n)}=\lim_{m \ra +\infty} (x^{(n+m)})^{p^m}$. This map then realizes an isomorphism
between $\mathcal{N}_K$ and $\e_K$, so that the two Galois groups $H_K$ and $G_{E_K}$ are
naturally isomorphic.

\begin{bibli}
For the theory of the field of norms in a much more general setting, see
Fontaine and Wintenberger's \cite{FW79} and Wintenberger's \cite{Wi83}. 
For the construction of the isomorphism
$\mathcal{N}_K \ra \e_K$ and its relation to Coleman series, see
Fontaine's appendix to \cite{BP94} and
Cherbonnier-Colmez's \cite{CC99}.
\end{bibli}

\subsubsection{$(\phi,\Gamma)$-modules}\label{pgmod}
By combining the construction of \ref{charp} and the theory of the field of norms, we see
that we have an equivalence of categories:
\[ \{ \text{$p$-adic representations of $H_K$}\} \qquad
\longleftrightarrow
\qquad \{ \text{\'etale $\phi$-modules} \}. \]
We immediately deduce from this the equivalence
of categories we were looking for:
\[ \{ \text{$p$-adic representations of $G_K$}\} \qquad
\longleftrightarrow
\qquad \{ \text{\'etale $(\phi,\Gamma_K)$-modules} \}. \]
One associates to $V$ the \'etale $\phi$-module $\dfont(V) = (\mathbf{B} \otimes_{\Qp}
V)^{H_K}$, which is an \'etale $\phi$-module
endowed with the residual action of $\Gamma_K$:
it is a \emph{$(\phi,\Gamma_K)$-module}.
The inverse functor
is then given by $D \mapsto (\mathbf{B} \otimes_{\mathbf{B}_K} D)^{\phi=1}$.

In general, it is rather hard to write down the $(\phi,\Gamma)$-module associated to a
representation $V$. We can therefore only give a few trivial examples, such as 
$\dfont(\Qp(r))=\mathbf{B}_F(r)$. See also the examples in
\ref{ordipg} and at the end of \ref{crys}. 

\begin{bibli}
The original theory of $(\phi,\Gamma)$-modules is the subject of Fontaine's 
\cite{Fo91}. It has been modified a bit 
by Cherbonnier and Colmez in \cite{CC99}, whose constructions we have followed.
\end{bibli}

\subsubsection{A flight of fancy}\label{crazy}
The above equivalence of categories
shows that in order to describe a $p$-adic representation of $G_K$, one only needs to
write down two matrices in $\on{GL}(d,\mathbf{A}_K)$: the matrix $P$ of $\phi$ and the matrix
$G$ of a topological generator $\gamma$ of $\Gamma_K$ (which is procyclic except
in some cases when $p=2$). 
These two matrices must satisfy the relation $\gamma(P)G=\phi(G)P$.
Let us elaborate a bit on this. In order to describe a 
$p$-adic representation of $H_K$, we only need to
write down the matrix $P \in \on{GL}(d,\mathbf{A}_K)$ of $\phi$. It is well defined up to
changing basis, that is if we let $\on{GL}(d,\mathbf{A}_K)$ act on itself by $X \cdot M =
\phi(X)MX^{-1}$, then the class of $P$ in $\on{GL}(d,\mathbf{A}_K) /
\sim$ (two matrices $M_1$ and $M_2$ being equivalent if
$M_1 = \phi(X)M_2 X^{-1}$)
is well defined. The latter coset space is endowed with an
action of $\Gamma_K$, which sends $\overline{M}$ to $\gamma(\overline{M})$, and a
representation of
$H_K$, given by a matrix $P$, extends to $G_K$ if and only if the image of $P$ in
$\on{GL}(d,\mathbf{A}_K) / \sim$ lies in
$(\on{GL}(d,\mathbf{A}_K) / \sim)^{\Gamma_K}$. In that case, there
are many different ways of extending that action (one needs to choose a lift of the class of
$P$). Of course, this is just an explicit way of looking at the inflation-restriction
sequence (which is an exact sequence of pointed sets):
\begin{multline*} \{1\} \ra H^1(\Gamma_K,\on{GL}(d,\mathbf{A}_K)) \ra
H^1(G_K,\on{GL}(d,\mathbf{A}_K)) \\ \ra
H^1(H_K,\on{GL}(d,\mathbf{A}_K))^{\Gamma_K} \ra \{1\}. \end{multline*}

Here is a possible application. Given two representations of $G_K$ which are close enough,
can we find a ``continuous path'' linking the two? One can easily find a continuous path of
representations of $H_K$, by interpolating the matrices of $\phi$. If we could impose a
minimality condition on that path (such as: being a geodesic), then the path would be
unique, and because its endpoints are in $(\on{GL}(d,\mathbf{A}_K) /
\sim)^{\Gamma_K}$, so would all its intermediate points: that is,
there would be a path of representations of $G_K$. The same approach should work equally
well with crystalline representations (see \ref{crys}) because they are described by
$(\on{GL}(d,\mathbf{A}_F^+) / 
\sim )^{\Gamma_F}$ where 
$M_1 \sim M_2$ if $M_1 = \phi(X)M_2 X^{-1}$
with $X \in \on{Id}+\pi \on{M}(d,\mathbf{A}_F^+)$.

\subsubsection{Computation of Galois cohomology}\label{galcoh}
Since the category of \'etale (i.e. slope $0$) 
$(\phi,\Gamma)$-modules is equivalent to that of $p$-adic
representations, it should be possible to recover all properties of $p$-adic
representations in terms of $(\phi,\Gamma)$-modules. For example, Herr showed in his thesis
how one could compute the Galois cohomology of $V$ from $\dfont(V)$. Let $H^i(K,V)$ denote
the groups of continuous cohomology of $V$. Herr's main result is that one can recover the
$H^i(K,V)$'s from $\dfont(V)$.

Let $\Delta_K$ be the torsion subgroup of
$\Gamma_K$; since $\Gamma_K$ is an open subgroup of $\Zp^*$,
$\Delta_K$ is a finite subgroup whose order divides $p-1$ (or $2$ if $p$=2). Let 
$p_{\Delta}$ be the idempotent defined by $p_{\Delta}=\frac{1}{|\Delta_K|}
\sum_{\delta\in\Delta_K} \delta$ so that if 
$M$ is a $\Zp[[\Gamma_K]]$-module, then
$p_{\Delta}$ is a projection map from $M$ to $M^{\Delta_K}$
(at least if $p \neq 2$). 
Let $\gamma$ be a topological
generator of $\Gamma_K/\Delta_K$.

Let $\dfontp(V)=\dfont(V)^{\Delta_K}$.
If $\alpha$ is a map $\alpha: \dfontp(V) \ra \dfontp(V)$ which commutes with $\Gamma_K$,
let $C_{\alpha,\gamma}(K,V)$ be the following complex :
\[ 0 \ra \dfontp(V) \overset{f}{\ra}
\dfontp(V)\oplus \dfontp(V)  \overset{g}{\ra}
\dfontp(V) \ra 0 \]
where $f(x)=((\alpha-1)x,(\gamma-1)x)$ and $g(x,y)=(\gamma-1)x-(\alpha-1)y$.

The cohomology of the complex $C_{\phi,\gamma}(K,V)$ is then naturally isomorphic to
the Galois cohomology of $V$. For example, we see immediately 
that $H^i(K,V)=0$ if $i\geq 3$. 

\begin{bibli}
This was proved by Herr in \cite{He00}. For various applications, see 
Herr's \cite{He00,He01,He02}, Benois'
\cite{Be00}, or \cite[chap VI]{LB00}, Cherbonnier-Colmez's \cite{CC99},
and Colmez's \cite{Co99}.
\end{bibli}

\Subsection{Overconvergent $(\phi,\Gamma)$-modules}\label{overcv}
Since the theory of $(\phi,\Gamma)$-modules is so good at dealing with $p$-adic
representations, we would like to be able to recover from $\dfont(V)$ the invariants
associated to $V$ by $p$-adic Hodge theory. This is the subject of the next chapter, on
reciprocity laws, but in this paragraph we will introduce the main technical tool, the ring
of overconvergent elements.

By construction, the field $\mathbf{B}$ is a subfield of 
\[ \bt = W(\et)[1/p] = \{
\sum_{k \gg -\infty} p^k [x_k],\ x_k \in \et \}. \] 
Let $\bdag{,r_n}$ be the subring of $\mathbf{B}$ defined as follows:
\[ \bdag{,r_n} = \{ x \in \mathbf{B},\ x= \sum_{k \gg -\infty} p^k [x_k],\ 
k+p^{-n}v_E(x_k) \ra +\infty \}. \]
The latter condition is equivalent to requiring that 
$\sum_{k \gg -\infty} p^k x_k^{(n)}$ converge in $\Cp$, which in turn is equivalent to
requiring that
$\sum_{k \gg -\infty} p^k [x_k^{p^{-n}}]$ converge in $\bdR^+$. 
If $e_K=[K_{\infty}:F_{\infty}]$, then one can show that 
there exists $\pi_K \in \mathbf{A}_K$ 
(see the end of \ref{charp})
such that for $n\gg 0$, the
invariants of $\bdag{,r_n}$ under the action of $H_K$ are given by
\begin{multline*}
(\bdag{,r_n})^{H_K} = \bdag{,r_n}_K =
\{\sum_{k=-\infty}^{+\infty} a_k \pi_K^k,\ 
\text{where $\sum_{k=-\infty}^{+\infty} a_k X^k$} 
\\ \text{is convergent and bounded on $p^{-1/e_K
r_n} \leq |X| < 1$} \}. \end{multline*} 
If $K=F$ (so that $e_K=1$), then one can take $\pi_F=\pi$, and 
the above description is valid for all $n \geq 1$.

A $p$-adic representation is said to be \emph{overconvergent} if, 
for some $n \gg 0$, $\dfont(V)$
has a basis consisting of elements of 
$\ddag{,r_n}(V)=(\bdag{,r_n} \otimes_{\Qp} V)^{H_K}$. This is equivalent to
requiring that there exist a basis of $\dfont(V)$ in which
$\on{Mat}(\phi) \in \on{M}(d, \bdag{,r_n}_K)$ for some $n \gg 0$.  

The main result on $p$-adic representations (or, equivalently, on \'etale
$(\phi,\Gamma)$-modules) is a theorem of Colmez and Cherbonnier which shows that every
$p$-adic representation of $G_K$ 
(equivalently, every \'etale $(\phi,\Gamma)$-module) is overconvergent. It is
not true that every \'etale $\phi$-module is overconvergent, and their proof uses the action
of $\Gamma$ in a crucial way. For instance, there is no such result in
the characteristic $p$ theory.

\begin{bibli}
The above result is the main theorem of Cherbonnier-Colmez's \cite{CC98}. 
Most applications of
$(\phi,\Gamma)$-modules to $p$-adic Hodge theory make use of it. If $V$ is absolutely
crystalline, then one can say more about the periods of $\dfont(V)$, see 
Colmez's \cite{Co99} and \cite[3.3]{LB01}. 
\end{bibli}

\[ \heartsuit\clubsuit\diamondsuit\spadesuit \]

\newpage

\section{Reciprocity laws for $p$-adic representations}
\Subsection{Overview}
\subsubsection{Reciprocity laws in class field theory}\label{recip}
The aim of this chapter is to give constructions relating the theory of
$(\phi,\Gamma)$-modules to $p$-adic Hodge theory. The first thing we'll do is explain why we
(and others) have chosen to call such constructions reciprocity laws. Recall that, in its
simplest form, the aim of class field theory is to provide a description of
$\on{Gal}(K^{\mathrm{ab}}/K)$, 
where $K$ is a field. For example, if $K$ is a local field, then one
has for every finite extension $L/K$ the norm residue symbol $(\cdot,L/K):K^* \ra
\on{Gal}(L/K)^{\mathrm{ab}}$, 
which is a surjective map whose kernel is $N_{L/K}(L^*)$. This is a
form of the local reciprocity law, and the aim of explicit reciprocity laws is to describe
(explicitly!) the map $(\cdot,L/K)$ (more precisely, the Hilbert symbol). 
For example, a theorem of Dwork shows that if $\zeta$ is a $p^n$-th root of unity,
then one has $(u^{-1},\Qp(\zeta)/\Qp)\cdot \zeta = \zeta^u$.

Let $V=\Qp(1)$, which is the Tate module of the multiplicative group $\mathbf{G}_m$.
The classical reciprocity map relates the tangent space $\ddR(V)$ of $\mathbf{G}_m$ to the
Galois cohomology $H^1(G_K,V)$. This is why we call a \emph{reciprocity map} those maps which
relate Galois cohomology and $p$-adic Hodge theory. Since the Galois cohomology 
of $V$ naturally 
occurs in $(\phi,\Gamma)$-modules, it is natural to call 
``reciprocity map'' those maps which relate $(\phi,\Gamma)$-modules and $p$-adic Hodge
theory.

This is the aim of this chapter: we will show how to recover $\dcris(V)$ or $\dst(V)$ from
$\dfont(V)$ and how to characterize de Rham representations. As an application, we will
prove Fontaine's monodromy conjecture.

\begin{bibli}
The first important constructions relating $(\phi,\Gamma)$-modules and $p$-adic Hodge theory
were carried out in Cherbonnier-Colmez's
\cite{CC99}, and are closely related to Perrin-Riou's exponential, as in
her \cite{BP94} and Colmez's \cite{Co98}. See also \cite{LBbk} for
explicit construction of Bloch-Kato's maps.
\end{bibli}

\Subsection{A differential operator on $(\phi,\Gamma)$-modules}
In order to further relate the theory of $(\phi,\Gamma)$-modules to $p$-adic Hodge theory,
we will need to look at the action of the Lie algebra of $\Gamma$ on $\ddag{}(V)$. On,
$\bdag{}_K$, it acts through a differential operator $\nabla$, given by
$\nabla=\log(\gamma)/\log_p(\chi(\gamma))$, which acts by
$\nabla(f(\pi))=\log(1+\pi)(1+\pi)df/d\pi$. We see that $\nabla(f(\pi)) \notin \bdag{}_K$,
and so it is necessary to extend the scalars to
\begin{multline*} 
\bnrig{,r_n}{,K} = \{ f(\pi_K) = \sum_{k=-\infty}^{+\infty} a_k \pi_K^k,\ \\
\text{where $f(X)$ is convergent on $p^{-1/e_K r_n} \leq |X| < 1$} \}.
\end{multline*}

The definition is almost the same as that of $\bdag{,r_n}_K$, but we have dropped the
boundedness condition. 
A typical element of $\bnrig{,r_n}{,K}$ is $t=\log(1+\pi)$.
We see that $\bnrig{,r_n}{,K}$ is a Fr\'echet space, with all the
norms given by the sup norms on ``closed'' annuli, and that it contains $\bdag{,r_n}_K$ as a
dense subspace. The union $\bnrig{}{,K} = \cup_{n=0}^{+\infty} \bnrig{,r_n}{,K}$ 
is the \emph{Robba ring} $\mathcal{R}_K$ of $p$-adic differential equations, and
$\mathcal{E}^{\dagger}_K = \bdag{}_K$ is the subring of $\mathcal{R}_K$ consisting of those
functions which are bounded. The $p$-adic completion of $\mathcal{E}^{\dagger}_K =
\bdag{}_K$ is $\mathcal{E}_K = \mathbf{B}_K$.

This being done, we see that the formula $\nabla_V=\log(\gamma)/\log_p(\chi(\gamma))$ gives
the action of $\on{Lie}(\Gamma_K)$ on $\dnrig{}(V) = \bnrig{}{,K} \otimes_{\bdag{}_K}
\ddag{}(V)$. Unfortunately, the action of $\on{Lie}(\Gamma_K)$ on $\bnrig{}{,K}$ is not very
nice, because $\nabla(f(\pi))=\log(1+\pi)(1+\pi)df/d\pi$ and this operator has zeroes at all
the $\zeta-1$ with $\zeta \in \mu_{p^\infty}$. 
In particular, it is not a basis of $\Omega^1_{\bnrig{}{,K}}$ and it is
not the kind of differential operator that fits within the framework of $p$-adic
differential equations. The ``right'' differential operator is $\partial_V =
\frac{1}{\log(1+\pi)} \nabla_V$, but this operator acting on $\dnrig{}(V)$ has poles at all
the $\zeta-1$. In the following paragraphs, we will see that
one can ``remove'' these poles exactly when $V$ is de Rham.

\begin{bibli}
See \cite{LB00,LB01} for detailed constructions and the basic
properties of those rings and operators.
\end{bibli}

\Subsection{Crystalline and semi-stable representations}
\subsubsection{Construction of $\dcris(V)$ and of $\dst(V)$}\label{dcrisst}
We will start by studying the action of $\Gamma_K$ on $\dnrig{}(V)$, and our main result is
that $\dcris(V) = (\dnrig{}(V)[1/t])^{\Gamma_K}$, in a sense which will be made precise
below. In addition, one can define 
$\bnst{}{}=\bnrig{}{}[\log(\pi)]$ with the obvious actions of $\phi$ and $\Gamma_K$, and we
shall also see that $\dst(V) = (\dnst{}(V)[1/t])^{\Gamma_K}$. If the Hodge-Tate weights of
$V$ are negative (if $V$ is \emph{positive}), then $\dcris(V) = \dnrig{}(V)^{\Gamma_K}$
and $\dst(V) = \dnst{}(V)^{\Gamma_K}$.

Recall that $\bst$ is a subring of $\bdR$ equipped with a Frobenius. The periods of $V$
are the elements of $\bst$ which ``occur'' in the coefficients of $\dst(V)$, they form a
finite dimensional $F$-vector subspace of $\bst$, stable by Frobenius. Therefore, these
periods live in $\cap_{n=0}^{+\infty} \phi^n(\bst^+)[1/t]$. 

The main strategy for comparing the theory of $(\phi,\Gamma)$-modules and $p$-adic Hodge
theory is to construct a rather large ring $\btrig{}{}$, 
which contains $\bdag{}$, $\bnrig{}{,K}$ and
$\cap_{n=0}^{+\infty} \phi^n(\bcris^+)$ so that
$\bnrig{}{,K} \otimes_{\bdag{}_K} \ddag{}(V) \subset \btrig{}{} \otimes_{\Qp} V$
and $\dcris(V) \subset (\btrig{}{} \otimes_{\Qp} V)^{G_K}$. The result alluded to above, for
positive crystalline representations, is that the two $F$-vector subspaces of
$\btrig{}{} \otimes_{\Qp} V$, 
$\dcris(V)$ and $\dnrig{}(V)^{\Gamma_K}$, actually coincide. 
This means that if $V$ is
crystalline, then the Frobenius $\phi$ on $\ddag{}(V)$ has a rather special form.
We'll give an informal justification for the above result in the next paragraph.

\subsubsection{Rings of periods and limits of algebraic functions}\label{algfun}
First of all, one should think
of most rings of periods as rings of ``limits of algebraic functions'' on certain subsets
of $\Cp$.
For example, the formula $\BB=\hat{\BB}_F^{\mathrm{unr}}$ 
tells us that $\BB$ is the ring of limits
of (separable) algebraic functions on the boundary of the open unit disk. The ring $\bt$ is
then the ring of all limits of algebraic functions on the boundary of the open unit disk.

Heuristically, one should view other rings in the same fashion: the ring $\bcris^+$ ``is''
the ring of limits of algebraic functions on the disk $D(0,|\eps^{(1)}-1|_p)$, and
$\bmax^+$ ``is'' the ring of limits of algebraic functions on a slightly smaller disk
$D(0,r)$. One should therefore think of $\phi^n(\bcris^+)$ as the ring of limits of algebraic
functions on the disk $D(0,|\eps^{(n)}-1|_p)$, and finally $\btrigplus{}$ ``is''  the ring
of limits of algebraic functions on the open unit disk $D(0,1)$.

Similarly, $\btrig{,r}{}$ ``is'' the ring of limits of algebraic functions on an annulus
$C[s,1[$, where $s$ depends on $r$, and $\phi^{-n}(\btrig{,r}{})$ 
``is'' the ring of limits of algebraic functions on an annulus $C[s_n,1[$, where $s_n \ra
0$, so that $\cap_{n=0}^{+\infty} \phi^{-n}(\btrig{,r}{})$ ``is'' the ring
of limits of algebraic functions on the open unit disk $D(0,1)$ minus the origin;
furthermore, if an element of that ring satisfies some simple growth properties near the
origin, then it ``extends'' to the origin (remember that in complex analysis, a holomorphic
function on $D(0,1^-)-\{0\}$ which is bounded near $0$ extends to a  holomorphic
function on $D(0,1^-)$). 

As for the ring $\bdR^+$, it behaves like a ring of local functions around a circle (in
particular, there is no Frobenius map defined on it). Via the map $\phi^{-n} : \bnrig{,r_n}{}
\ra \bdR^+$, we have for $n \geq 1$
a filtration on $\bnrig{,r_n}{}$, which corresponds to the order of
vanishing at $\eps^{(n)}-1$. For instance, we can now give a short solution to the exercise
in paragraph \ref{exseq}: 
given a sequence $r_n$ of integers, let $q=\phi(\pi)/\pi$ and set $x_r =
\pi^{r_0} \prod_{n=1}^{+\infty} \phi^{n-1}(q/p)^{r_n}$. This infinite product converges to a
``function'' whose order of vanishing at $\eps^{(n)}-1$ is exactly $r_n$.

\subsubsection{Regularization and decompletion}\label{diffst}
We shall now justify the above results on $\dcris(V)$. The analogous results on $\dst(V)$
follow by adding $\log(\pi)$ everywhere. We've already seen that the periods of
positive crystalline representations live in $\btrigplus{}$ (if we don't assume that $V$ is
positive, then they live in $\btrigplus{}[1/t]$).

The elements of $(\btrig{}{} \otimes_{\Qp} V)^{G_K}$ form a finite dimensional $F$-vector
space, so that there is an $r$ such that $(\btrig{}{} \otimes_{\Qp} V)^{G_K} = (\btrig{,r}{}
\otimes_{\Qp} V)^{G_K}$, and furthermore this $F$-vector space is stable by Frobenius, so
that the periods of $V$ (in this setting) not only live in $\btrig{,r}{}$ but actually in 
$\cap_{n=0}^{+\infty} \phi^{-n}(\btrig{,r}{})$ and they also satisfy some simple growth
conditions (depending, say, on the size of $\det(\phi)$), which ensure that they too can be
seen as limits of algebraic functions on the open unit disk $D(0,1)$, that is as
elements of $\btrigplus{}$. In particular, we have  
$(\btrig{}{} \otimes_{\Qp} V)^{G_K} = (\btrigplus{} \otimes_{\Qp} V)^{G_K}$.
This is what we get by \emph{regularization} (of the periods).

It's easy to show that $(\btrig{}{} \otimes_{\Qp} V)^{H_K} = \btrig{}{,K}
\otimes_{\bnrig{}{,K}} \dnrig{}(V)$, and 
the last step is to show that 
$(\btrig{}{,K} \otimes_{\bnrig{}{,K}} \dnrig{}(V))^{G_K} = \dnrig{}(V)^{G_K}$.
This is akin to a \emph{decompletion} process, 
going from $\btrig{}{,K}$ to $\bnrig{}{,K}$. The
ring extension 
$\btrig{}{,K} / \bnrig{}{,K}$
looks very much like $\hat{K}_{\infty}/K$, so that
by using Colmez's decompletion maps, 
which are analogous to Tate's $\on{pr}_{K_n}$ maps from paragraph \ref{ast},
one can finally show that in fact,
$\dcris(V) = (\bnrig{}{,K} \otimes_{\bdag{}_K} \ddag{}(V))^{G_K}$. In particular, $V$ is
crystalline if and only if $(\bnrig{}{,K} \otimes_{\bdag{}_K} \ddag{}(V))^{G_K}$ is a 
$d$-dimensional $F$-vector space.

\begin{bibli}
See \cite{LB00,LB01}. For decompletion maps and the ``Tate-Sen'' conditions, see Colmez's
explanations in \cite{Co01}.
\end{bibli}

\Subsection{De Rham representations}
In the previous paragraph, we have shown how to recognize crystalline and semi-stable
representations in terms of their $(\phi,\Gamma)$-modules. We shall now do the same for de
Rham representations, and show that a representation $V$ is positive de Rham if and only if 
there exists a free $\bnrig{}{,K}$-submodule of rank $d$ of $\dnrig{}(V)$, 
called $\ndr(V)$, which
is stable by the operator $\partial_V$
(when $V$ is not positive, then
$\ndr(V) \subset \dnrig{}(V)[1/t]$).
Of course, when $V$ is 
crystalline or semi-stable,
one can simply take $\ndr(V)=\bnrig{}{,K} \otimes_F \dcris(V)$ or
$\ndr(V)=(\bnst{}{,K} \otimes_F \dst(V))^{N=0}$. 

\subsubsection{Construction of $\ndr(V)$}\label{ndr}
In general, let us give an idea of how one can construct $\ndr(V)$. In the paragraph
\ref{bdrsen}, we recalled Fontaine's construction of ``Sen's theory for $\bdR^+$''. The map
$\phi^{-n}$ sends
$\dnrig{,r_n}(V)$ into $(\bdR^+ \otimes_{\Qp} V)^{H_K}$, which should be thought of as
``localizing at $\eps^{(n)}-1$'' in geometrical terms. The module $\ddif^+(V)$ of Fontaine
is then equal to $K_{\infty}[[t]] \otimes_{\phi^{-n}(\bnrig{,r_n}{,K})} \dnrig{,r_n}(V)$.
Recall that Fontaine has shown that a positive $V$ is de Rham and if and only if the
connexion $\nabla_V$ has a full set of sections on $\ddif^+(V)$ (in which case the kernel
of the connexion is $K_{\infty} \otimes_K \ddR(V)$).  In geometrical terms, this
means that if $V$ is positive and de Rham, then $\nabla_V$ has some ``local'' solutions
around the $\eps^{(n)}-1$. In that case, one can glue all of those solutions together to
obtain $\ndr(V)$. More precisely, there exists $n_0 \gg 0$ and
$r\gg 0$ such that we have $\ndr(V) =
\bnrig{}{,K} \otimes_{\bnrig{,r}{,K}} N_r(V)$ where $N_r(V)$  is the set of $x \in
\dnrig{,r}(V)$ such that for every $n \geq n_0$, one has $\phi^{-n}(x) \in
K_n[[t]] \otimes_K \ddR(V)$. It's easy to see that $N_r(V)[1/t] = \dnrig{,r}(V)[1/t]$ and
that $N_r(V)$ is a closed (for the Fr\'echet topology) $\bnrig{,r}{,K}$-submodule of
$\dnrig{,r}(V)$. The fact that $N_r(V)$ is free of rank $d$ then follows form the following
fact: if $M \subset (\bnrig{,r}{,K})^d$ is a closed submodule, such that
$\on{Frac} \bnrig{,r}{,K} \otimes_{\bnrig{,r}{,K}}  M = (\on{Frac} \bnrig{,r}{,K})^d$, then
$M$ is free of rank $d$.

One can then show that $\ndr(V)$ is uniquely determined by the requirement that it be free
of rank $d$ and stable by $\partial_V$, so that in particular $\phi^* \ndr(V) = \ndr(V)$.

We therefore have the following theorem: if $V$ is a de Rham representation, then there
exists $\ndr(V) \subset 
\dnrig{}(V)$, a $\bnrig{}{,K}$-module
free of rank $d$, stable by
$\partial_V$ and $\phi$, such that $\phi^* \ndr(V) = \ndr(V)$. Such an object is by
definition a $p$-adic differential equation with Frobenius structure
(see \ref{kam} below).

Using this theorem, one can construct a faithful and essentially surjective 
exact $\otimes$-functor from the category of de Rham representations to the category of
$p$-adic differential equations with a Frobenius structure.

\begin{bibli}
The above theorem is the main result of \cite{LB01}. The result on
closed submodules of $(\bnrig{,r}{,K})^d$ is proved in \cite[4.2]{LB01}, see also
\cite{FO67}. 
\end{bibli}

\subsubsection{Example: $\Cp$-admissible representations}\label{tsuzuki}
Let us give an example for which it is easy to characterize $\ndr(V)$. We've already seen that
when $V$ is crystalline or semi-stable,
one can take $\ndr(V)=\bnrig{}{,K} \otimes_F \dcris(V)$ or
$\ndr(V)=(\bnst{}{,K} \otimes_F \dst(V))^{N=0}$. Another easy case is when $V$ is
$\Cp$-admissible. This was one of the examples in \ref{strat} where we mentioned Sen's
result: a representation $V$ is $\Cp$-admissible if and only if it is potentially
unramified. We'll give a proof of that result which relies on a theorem of Tsuzuki on
differential equations.

Let $V$ be a $\Cp$-admissible representation. This means that $\Cp \otimes_{\Qp} V =
\Cp \otimes_K (\Cp \otimes_{\Qp} V)^{G_K}$, so that $V$ is Hodge-Tate and all its weights
are $0$. In particular, Sen's map $\Theta_V$ is zero. Since we recovered Sen's map 
from $\nabla_V$ by localizing at $\eps^{(n)}-1$, this implies that the coefficients of a
matrix of $\nabla_V$ are holomorphic functions which are $0$ at  $\eps^{(n)}-1$ for all $n
\gg 0$. These functions are therefore multiples of $t=\log(1+\pi)$ in $\bnrig{}{,K}$ and so 
$\nabla_V (\dnrig{}(V)) \subset \log(1+\pi) \dnrig{}(V)$ so that one can take
$\ndr(V)=\dnrig{}(V)$. 

The $\mathcal{R}_K$-module $\ndr(V)$ is then endowed with a differential operator
$\partial_V$ and a unit-root Frobenius map $\phi$ which is overconvergent. 
One can show that if $\phi$ is overconvergent, then so is $\partial_V$ (because $\phi$
regularizes functions). The module $\ndr(V)$ is therefore an \emph{overconvergent unit-root
isocrystal}, and Tsuzuki proved that these are potentially trivial (that is, they become
trivial after extending the scalars to $\mathcal{R}_L/\mathcal{R}_K$ for a finite extension
$L/K$). This implies easily enough
that the restriction of $V$ to $I_K$ is potentially trivial.

\begin{bibli}
See \cite[5.6]{LB01}. For Tsuzuki's theorem, see his \cite{TS99} and Christol's
\cite{GC00}. Sen's theorem was first proved in Sen's \cite{Sn73}.
\end{bibli}

\Subsection{The monodromy theorem}
\subsubsection{$\ell$-adic monodromy and $p$-adic monodromy}
\label{gro}
As was pointed out in the introduction, $\ell$-adic representations are forced to be
well-behaved, while the group $G_K$ has far too many $p$-adic representations.  
Over the years it became apparent that the only representations related to arithmetic
geometry were the de Rham representations (see \ref{tsuji}). 

In particular it was conjectured (and later
proved) that all representations coming from geometry were de Rham. Among these, some are
more pleasant, they are the semi-stable ones, which are the analogue of the $\ell$-adic
unipotent representations. Grothendieck has shown that all $\ell$-adic representations
are quasi-unipotent, and after looking at many examples, Fontaine was led to conjecture the
following  $p$-adic analogue of Grothendieck's $\ell$-adic monodromy theorem:
every de Rham representation is potentially semi-stable. We shall now give a proof of that
statement.

\begin{bibli}
An excellent reference throughout this section is Colmez's Bourbaki talk \cite{Co01}.
\end{bibli}

\subsubsection{$p$-adic differential equations}\label{kam}
A \emph{$p$-adic differential equation} 
is a module $M$, free of finite rank over the Robba ring
$\mathcal{R}_K$, equipped with a connexion $\partial_M : M \ra M$. We say that $M$ has a
\emph{Frobenius structure} 
if there is a semi-linear Frobenius $\phi_M : M \ra M$ which commutes with
$\partial_M$. A $p$-adic differential equation is said to be
\emph{quasi-unipotent} if there exists a finite extension 
$L/K$ such that $\partial_M$ has a full set of solutions on 
$\mathcal{R}_L[\log(\pi)] \otimes_{\mathcal{R}_K} M$. Christol and
Mebkhout extensively studied $p$-adic differential equations. 
Crew and Tsuzuki conjectured that every $p$-adic differential equation with a Frobenius
structure is quasi-unipotent.
Three independent proofs were given in the summer of $2001$.
One by Andr\'e, using Christol-Mebkhout's results and a
clever Tannakian argument. One by Kedlaya, who proved a
``Diedudonn\'e-Manin'' theorem for $\phi$-modules over
$\mathcal{R}_K$. And one by Mebkhout, relying on 
Christol-Mebkhout's results.

\begin{bibli}
We refer the reader to Christol and Mebkhout's surveys
\cite{CMpde,CM00} and Colmez's Bourbaki talk \cite{Co01} for 
enlightening discussions of $p$-adic differential equations.
The above theorem is proved independently in Andr\'e's \cite{An01}, 
Mebkhout's \cite{Mk01} and Kedlaya's
\cite{KK00}. See also Andr\'e's \cite{An00} for a beautiful discussion 
of a special case.
\end{bibli}

\subsubsection{The monodromy theorem}\label{mono}
Using the previous results, one can give a proof of Fontaine's monodromy conjecture. Let $V$
be a de Rham representation, then one can associate to $V$ a $p$-adic
differential equation $\ndr(V)$. By Andr\'e, Kedlaya, and Mebkhout's theorem, this
differential equation is quasi-unipotent. Therefore, there exists a finite extension $L/K$ 
such that $(\mathcal{R}_L[\log(\pi)] \otimes_{\mathcal{R}_K} \ndr(V) )^{G_L}$ is an
$F$-vector space of dimension $d$ and by the results of
paragraph \ref{diffst}, $V$ is potentially semi-stable. 

\begin{bibli}
See \cite[5.5]{LB01} for further discussion of the above result.
\end{bibli}

\subsubsection{Example: Tate's elliptic curve}\label{ordipg}
To finish this chapter, 
we will sketch this for Tate's elliptic curve (or indeed for all ordinary
elliptic curves). For simplicity,
assume that $k$ is algebraically closed. If
$q=q_0$ is the parameter associated to $E_q$, then there exists $q_n\in F_n = F(\eps^{(n)})$
such that
$N_{F_{n+1}/F_n}(q_{n+1})=q_n$
(this is the only place where we use the fact that $k$ is algebraically closed), 
and by a result of Coleman, there is a power series
$\on{Col}_q(\pi)$ such that $q_n = \on{Col}_q^{\sigma_F^{-n}}(\eps^{(n)}-1)$.
If $F_q(\pi)=(1+\pi)\on{dlog}\on{Col}_q(\pi)$, then $F_q(\pi) \in \pi^{-1} \OO_F[[\pi]]$ and
one can show that there is a basis $(a,b)$
of the $(\phi,\Gamma)$-module $\dfont(V)$ associated to $V$ such
that the action of $\Gamma_F=\langle \gamma \rangle$ is given by:
\[ \on{Mat}(\eta) =
\begin{pmatrix}
\chi(\eta) & \frac{1-\eta}{1-\gamma} F_q(\pi) \\
0 & 1 \end{pmatrix}\] 
Let $\nabla$ be the differential operator giving the action of the Lie algebra of $\Gamma_F$
on  power series, so that we have $(\nabla f)(\pi)=(1+\pi)\log(1+\pi)f'(\pi)$
(recall that $t=\log(1+\pi)$). 
The Lie algebra of
$\Gamma_F$ then acts on $D_{\mathrm{rig}}$ by an operator $\nabla_V$ given by
\[ \on{Mat}(\nabla_V)=\begin{pmatrix}
1 & \frac{\nabla}{1-\gamma} F_q(\pi) \\
0 & 0  \end{pmatrix}\] 
One then sees that $\partial_V(t^{-1}a) = 0$ and that $\partial_V(b)$ belongs to
$\bnrig{}{,F}(t^{-1}a)$, so that the $p$-adic differential equation $\langle t^{-1}a, b
\rangle$ is unipotent. This shows that $V$ is indeed semi-stable.

\begin{bibli}
The extensions of $\Qp$ by $\Qp(1)$ are important and also a source of explicit examples.
They are related to Kummer theory as in paragraph \ref{kummer}, and Coleman series as
above, among other topics. Some interesting computations can be found 
in Cherbonnier-Colmez's \cite[V]{CC99}.
\end{bibli}

\[ \spadesuit\heartsuit\clubsuit\diamondsuit \]

\newpage

\section{Applications and further topics}
In this chapter, we will mention a few recent applications of the
results which we explained above.

\Subsection{Limits of crystalline representations}
In this section, we will give an application of the above results to crystalline
representations. We'll assume throughout that $K=F$. The  result which
we will prove is that 
if $a,b \in \ZZ$ and $V_i = \Qp \otimes_{\Zp} T_i$ is a
sequence of crystalline representations of $G_F$ with Hodge-Tate weights in $[a;b]$, 
such that $T_i/p^i \simeq T_{i+1}/p^i$, and $V= \Qp \otimes_{\Zp}  T$ with $T = \lim_i
T_i/p^i$, then
$V$ is itself crystalline. Note that if $b-a
\leq p-1$, then the result trivially follows from Fontaine-Laffaille
\cite{FL82}. 
As an example of an application, Breuil has shown 
that some representations attached to Hilbert modular forms are
naturally limits of crystalline representations. 

The bound on the Hodge-Tate weights is important, for example if $s_i$ is a sequence of
integers such that $s_i = s_{i+1} \mod{p^{i-1}(p-1)}$, and if $s= \lim_{i\ra +\infty} s_i 
\in \projlim_i \ZZ/p^{i-1}(p-1)$ but not in
$\ZZ$, then $V_i=\Qp(s_i)$ is crystalline but $V=\Qp(s)$ is not
(because $s$ should be its Hodge-Tate weight, but is not an
integer). See also the remarks in \ref{famil}.

\subsubsection{Crystalline representations and finite height representations}\label{crys}
We can twist $V$ and the $V_i$ by $-b$, and we will therefore 
assume throughout this paragraph that 
the Hodge-Tate weights of $V$ are $\leq 0$.

If $V$ is a crystalline representation, 
we can say a lot about the $(\phi,\Gamma)$-module attached
to $V$. Let $\aplus_F=\OO_F[[\pi]]$ and $\bplus_F=\OO_F[[\pi]][1/p]$.
Using the fact that $\dcris(V)  = \dnrig{}(V)^{\Gamma_F}$,
and a result of Kedlaya on factorization of matrices in
$\bnrig{}{,F}$, one can show that if $V$ is a crystalline 
representation of $G_F$, 
and $T \subset V$ is a lattice of $V$, then
there exists a free of rank $d$
$\aplus_F$-module $N(T) \subset \ddag{}(V)$, 
stable by $\phi$ and $\Gamma_F$, such that 
$\Gamma_F$ acts trivially on $N(T)/\pi$ and such that
$\ddag{}(T)= \adag{}_F \otimes_{\aplus_F} N(T)$ so that
$\ddag{}(V)= \bdag{}_F \otimes_{\aplus_F} N(T)$. 
Conversely, the existence of such a $N(T)$
implies that $V$ is crystalline.
In other words, $\ddag{}(V)$ has a
basis in which the matrices of $P$ of $\phi$ and $G$ of $\gamma$ are
in $\on{M}(d,\aplus_F)$ and $G = \on{Id} \mod{\pi}$. One can then show that
$F \otimes_{\OO_F} N(T)/\pi \simeq \dcris(V)$ as $\phi$-modules.

We can also recover $\dcris(V)$ from $N(T)$ in 
the above fashion. Let $\bnrigplus{,F}$ be the set of
power series $\sum_{k=0}^{+\infty} a_k X^k$ which are 
holomorphic on the open unit disk. We have $\dcris(V) =
(\bnrigplus{,F} \otimes_{\aplus_F} N(T))^{\Gamma_F}$.

Let us give an example of all this. Consider the $2$-dimensional
filtered $\phi$-module $D=\dcris(V)$ generated by 
$\overline{e}$, $\overline{f}$ such that $\on{Fil}^0 D=D$, 
$\on{Fil}^1 D = F \overline{e}$, 
$\on{Fil}^2 D=0$, with $\phi(\overline{e})=p\overline{f}$ 
and $\phi(\overline{f})=\overline{e}$. 
It comes from a supersingular elliptic
curve. The module $N(T)$ associated to $D$ is then 
the free $\aplus_F$-module generated by $e$ and $f$ with
$\phi(e)=\Phi_1(1+\pi) f$ and $\phi(f)=e$ where
$\Phi_n(1+\pi)=((1+\pi)^{p^n}-1)/((1+\pi)^{p^{n-1}}-1)$. 
The action of $\Gamma_F$ is given by 
\[ \gamma(e)=\frac{\log^-(1+\pi)}{\gamma(\log^-(1+\pi))} e \quad\text{and}\quad 
\gamma(f)=\frac{\log^+(1+\pi)}{\gamma(\log^+(1+\pi))}f \] where
\[ \log^+(1+\pi) =\prod_{n \geq 0} \frac{\Phi_{2n}(1+\pi)}{p},\quad 
\log^-(1+\pi) =\prod_{n \geq 0} \frac{\Phi_{2n+1}(1+\pi)}{p} \]
so that $t=\log(1+\pi)=\pi \log^-(1+\pi) \log^+(1+\pi)$. It is then easy to
check that one can recover $\overline{e}$, $\overline{f}$ either as
the images modulo $\pi$ of $e$, $f$ or as $\overline{e} = \log^-(1+\pi) e$ and
$\overline{f} = \log^+(1+\pi) f$ as elements of  $\dcris(V) =
(\bnrigplus{,F} \otimes_{\aplus_F} N(T))^{\Gamma_F}$.

\subsubsection{The determinant of $\phi$ on $N_T$}\label{detphi}
The next ingredient we need is to compute $\det(P)$ where $P$ is the
matrix of $\phi$. Since $\bnrigplus{,F} \otimes_F \dcris(V) \subset
\bnrigplus{,F} \otimes_{\aplus_F} N(T)$, and we know the determinant
of $\phi$ on $\dcris(V)$, it will be enough to compute the elementary
divisors $\lambda_1, \cdots, \lambda_d$ of  
$\bnrigplus{,F} \otimes_F \dcris(V) \subset
\bnrigplus{,F} \otimes_{\aplus_F} N(T)$. 
The ideal generated by an analytic function is characterized by the
zeroes of that function. The first observation is that
the ideals $\lambda_i$ are stable by $\Gamma_F$, and this shows that
the only zeroes of the $\lambda_i$ are the $\zeta-1$, where $\zeta \in
\mu_{p^\infty}$. We know that $\Gamma_F$ acts trivially on $N(T)/\pi$ so
that the $\lambda_i$ are $\neq 0$ at $0$. If $\zeta^{p^n}=1$ with $n
\geq 1$, then
localizing at $\zeta-1$ is equivalent to 
\[ \begin{CD}
\iota_n : \bnrigplus{,F} \otimes_{\aplus_F} N(T) 
@>{\phi^{-n}}>> (\bdR^+ \otimes_{\Qp}
V)^{H_F}  @>{\theta}>> (\Cp \otimes_{\Qp} V)^{H_F}. \end{CD} \] 

Since $V$ is Hodge-Tate, we know what the action of $\Gamma_F$ on
$(\Cp \otimes_{\Qp} V)^{H_F}=\oplus_{j=1}^d \hat{F}_{\infty}(r_j)$ 
looks like and we can 
show that for every $\gamma \in \Gamma_n = 
\on{Gal}(F_{\infty}/F_n)$, the matrix of $\gamma$ on
the $F_n$-vector space $F_n \otimes_{\OO_{F_n}} N(T) / \Phi_{n}(1+\pi)$ 
is diagonalizable, with eigenvalues
$\{\chi^{r_j}(\gamma)\}_{j=1\cdots d}$. 
This tells us that at $\zeta-1$, $\lambda_i$ has a zero of order
$-r_i$ and therefore that $\lambda_i=(\log(1+\pi)/\pi)^{-r_i}$ so that
$(\det(P))=(\Phi_1(1+\pi)^{-t_H(V)})$ where $t_H(V)=\sum_{j=1}^d r_j$.

To summarize, we know that a representation 
$V$ is positive crystalline if and only if
its $(\phi,\Gamma)$-module contains an 
$\aplus_F$-module $N(T)$ free of rank $d$, 
stable by $\phi$ and the action of $\Gamma_F$, such that
$\Gamma_F$ acts trivially on $N(T)/ \pi$  and such 
that the determinant of $\phi$ is the ideal generated by
$\Phi_1(1+\pi)^{-t_H(V)}$. This means that a lattice $T$ of $V$ is 
determined by two matrices $P$ (of $\phi$) and $G$ 
(of $\gamma$, a topological generator of $\Gamma_F$) 
satisfying some simple conditions:
\[ P,G \in \on{M}(d,\OO_F[[\pi]]),\quad \gamma(P)G=\phi(G)P,
\quad \det(P) \mid \Phi_1(1+\pi)^{-t_H(V)}, \quad \pi
\mid G-\on{Id}.\] It's not hard to see that the 
set of such matrices is compact if one bounds
$t_H(V)$, so that if $V_i$ is a sequence of 
crystalline representations of $G_F$
with bounded weights, then its limit
will also be crystalline, essentially because 
one can take $N(T) = \projlim_i N(T_{n(i)})/p^{n(i)}$.

More precisely, by being a little more careful,
it is possible to show that if $T_1$ and $T_2$ are two lattices in two
crystalline representations, with Hodge-Tate weights in $[a;b]$, such that $T_1/p^n =
T_2/p^n$, then
$N(T_1) / p^{n-c} = N(T_2) / p^{n-c}$ 
with $c=v_p((b-a)!)$
provided that $p^{n-1}(p-1) > b-a$. 
It is then clear that one can (and should) take  $N(T) = \projlim_i N(T_i) / p^i$.

\begin{bibli}
The construction of $N(T)$ is originally due to
Wach \cite{Wa96}. The reference for this whole section is \cite{LB02}. 
Some partial results were known, see
the introduction to \cite{LB02} for an overview. The same results are expected to hold for
limits of semi-stable representations of $G_K$ (Berger-Colmez: work in progress).
See also \ref{famil}.
\end{bibli}

\Subsection{Bloch-Kato's exponential}
\label{expbk}
Recall that in \ref{exseq}, we stated the fundamental exact sequence:
\[ 0 \ra \Qp \ra \bcris^{\phi=1} \ra \bdR/\bdR^+ \ra 0 \]
By tensoring this sequence with a de Rham representation $V$ 
and taking $G_K$-invariants, we get the beginning of a long exact sequence:
\begin{multline*} 
0 \ra V^{G_K} \ra \dcris(V)^{\phi=1} \ra 
((\bdR/\bdR^+)\otimes_{\Qp} V)^{G_K} \\ \ra H^1(K,V) \ra
H^1(K,\bmax^{\phi=1}\otimes_{\Qp} V). \end{multline*}
Let $H^1_e(K,V)=\ker(H^1(K,V) \ra H^1(K,\bmax^{\phi=1}\otimes_{\Qp}
V))$. The above exact sequence then gives us a surjective map
from $\ddR(V)$ to $H^1_e(K,V)$ called \emph{Bloch-Kato's exponential} 
and denoted by
$\exp_V$.  Note that since $V$ is de Rham, one has 
\[ ((\bdR/\bdR^+)\otimes_{\Qp} V)^{G_K}=\ddR(V)/ \on{Fil}^0 \ddR(V). \]

If $\mathbf{G}$ is a formal group, and $V$ is its Tate module, then $\on{Tan}(\mathbf{G})$,
the tangent space to $\mathbf{G}$,
is naturally identified with $\ddR(V)/ \on{Fil}^0 \ddR(V)$ and Bloch-Kato's exponential is
just the composition of the usual exponential $\exp : \on{Tan}(\mathbf{G}) \ra \Qp
\otimes_{\Zp} \mathbf{G}(K)$ with the Kummer map $\delta:  \Qp \otimes_{\Zp} \mathbf{G}(K)
\ra H^1(K,V)$, hence the name.

Bloch and Kato have computed the dimension of $H^1_e(K,V)$ as a $\Qp$-vector space in case
$k$ is finite, and used their computations to show that if $V$ is a semi-stable
representation and $r\gg 0$, then the exponential $\exp_{V(r)}: \ddR(V(r)) \ra H^1(K,V(r))$
is an isomorphism. By using the results of the previous section, and some Galois cohomology,
we can show that $\exp_V: \ddR(V) \ra H^1(K,V)$ is an isomorphism for any perfect residual
field $k$, whenever the Hodge-Tate weights of $V$ are $\geq 2$. 

First of all, observe that
$H^1(K,V)$ classifies the extensions of $\Qp$ by $V$. It is
a simple exercise to show that if the weights of
$V$ are $\geq 1$, then all such extensions are de Rham. By the monodromy theorem, they are
potentially semi-stable, and an extension of semi-stable representations which is
potentially semi-stable is actually semi-stable 
(this is again a simple exercise in Galois cohomology).
This means that the natural map
$H^1(K,V) \ra H^1(K,\bst \otimes_{\Qp} V)$ is zero. If the weights of $V$ are $\geq 2$, then
(by looking at the action of $N$ and bearing in mind that $N\phi=p\phi N$), we show
that $N=0$ so that
the periods of these extensions are actually in $\bcris$, and
one can modify these periods a bit so that they all lie in $\bcris^{\phi=1}$ 
instead of just $\bcris$.
This shows that $\exp_V$ is surjective. Finally, its kernel is
$\on{Fil}^0 \ddR(V) + \dcris(V)^{\phi=1}$ which is
$0$ if the weights of $V$ are $\geq 1$.

\begin{bibli}
The reference for this section is \cite[6.2]{LB01}, see also the end of
paragraph \ref{exseq} for references on Bloch-Kato's exponential. 
\end{bibli}

\Subsection{Further topics}
We will now briefly mention a few interesting topics, that we do not
have the space (or the knowledge) to comment on.

\subsubsection{Integral $p$-adic Hodge theory}
One topic that we will say nothing about is integral $p$-adic Hodge theory. We refer the
reader to the paper 
of Breuil \cite{Br00} for an illuminating survey and further references.

\subsubsection{Perrin-Riou's exponential}
For arithmetic applications, it is desirable to understand how Bloch-Kato's exponential maps
$\exp_{V(k)}$ vary when $k$ varies continuously in weight space. An answer to that problem
was given by Perrin-Riou in \cite{BP94}; see also the reference in the
next paragraph. 
It is possible to give a construction of Perrin-Riou's map using
$(\phi,\Gamma)$-modules and some of the results mentioned 
above in this survey, see \cite{LBbk}.

\subsubsection{Construction of $p$-adic $L$-functions}
Another important application of $p$-adic Hodge theory is the construction of $p$-adic
$L$-functions, from the data of an Euler system. This topic is 
intimately linked to the preceding one, and it is the subject of 
a Bourbaki talk \cite{CoBo} by Colmez, to which we refer the reader. 

\subsubsection{Families of $p$-adic representations}\label{famil}
A \emph{family} of $p$-adic representations is 
(roughly speaking)
a $p$-adic space $X$ such that
to every point $x \in X$ one can attach a $p$-adic representation
$V_x$, and such that $x \mapsto V_x$ is analytic. For example, if
$\mathcal{R}$ is a Tate algebra, then a free $\mathcal{R}$-module with
an action of $G_K$ is a family of representations.

We would like to know what the subspaces of $X$, made up of
representations satisfying certain conditions look like. For example if 
we set
$X_{\mathrm{cris}} = \{ x \in X, V_x\ \text{is crystalline}\}$, 
what does $X_{\mathrm{cris}}$
look like? A good answer to that question could allow one to restrict
possibilities among all deformations of a representation $V$ when
looking for modular ones.

Sen has shown that the map $x \mapsto \Theta_{V_x}$, which to $V$
attaches Sen's operator $\Theta_V$, is analytic. Therefore 
if $J$ is a finite set of integers, and $X_J = \{ x \in X, V_x\ \text{has
Hodge-Tate weights in}\ J\}$, then $X_J$ is an analytic subspace of $X$.

In general, it is not true that $X_{\mathrm{cris}}$ or 
$X_{\mathrm{st}}$ is an analytic subspace, or
even a closed subspace of $X$. For example, if $X$ is the set of all
extensions $0 \ra \Qp(\eta) \ra E \ra \Qp \ra 0$, where $\eta$ is in
weight space $W$, then there is a natural map $\on{wt}: X \ra W$ and
$X_{\mathrm{st}}$ is made up of trivial representations, of
$\on{wt}^{-1}(1,2,3,\cdots)$, and of a one dimensional subspace of
$\on{wt}^{-1}(0)$. This is a dense subspace of $X$.

The theorem proved in \ref{crys} says that if $K=F$, then
$X_{J,\mathrm{cris}}$ is a
closed subspace of $X$. More generally, one can show that for any $K$,
the spaces 
$X_{J,\mathrm{cris}}$, $X_{J,\mathrm{st}}$ and $X_{J,\mathrm{dR}}$ are all
analytic subspaces of $X$.

\begin{bibli}
All of this is work in progress with Colmez. Kisin \cite{Ki02} 
has interesting
results on families of $p$-adic representations.
\end{bibli}

\[ \diamondsuit\spadesuit\heartsuit\clubsuit \]

\newpage

\section{Appendix}

\Subsection{Diagram of the rings of periods}
The following diagram summarizes the relationships between 
the different rings of periods. The arrows
ending with 
$\xymatrix@1{&\ar@{->>}[r]&}$ 
are surjective, the dotted arrow
$\xymatrix@1{&\ar@{.>}[r]&}$ 
is an inductive limit of maps defined on subrings
($\iota_n: \btst{,r_n}{} \ra \bdR^+$), 
and all the other ones are injective.
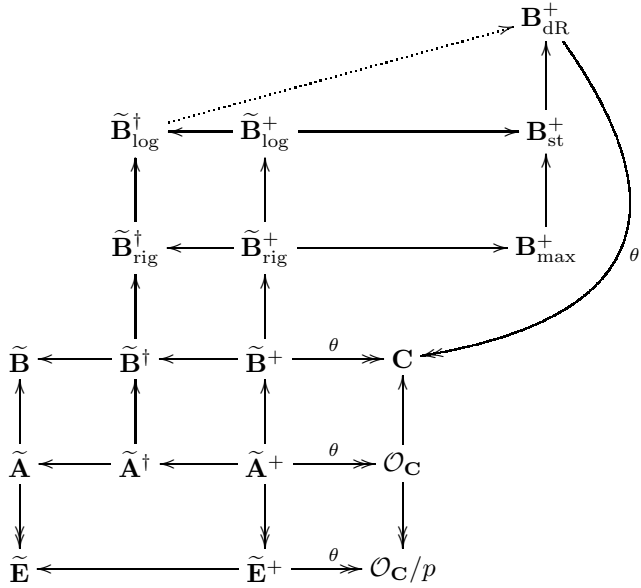
\[ \xymatrix{
& &  & & \bdR^+ \ar@{->>}@/^5pc/[lddd]^{\theta} & \\
& \btst{}{} \ar@{.>}[urrr] & \btstplus{} \ar[rr] 
\ar[l] & & \bst^+ \ar[u]  & \\
& \btrig{}{} \ar[u] & \btrigplus{} \ar[u] 
\ar[l] \ar[rr] & & \bmax^+ \ar[u]  & \\  
\bt & \btdag{} \ar[l] \ar[u] & \btplus \ar[l]
\ar[u] \ar@{->>}[r]^{\theta} & \Cp & & \\ 
\at \ar[u] \ar@{->>}[d] & \atdag{} \ar[l] \ar[u] 
& \atplus \ar[l] \ar[u] 
\ar@{->>}[d] \ar@{->>}[r]^{\theta} & 
\OO_{\Cp} \ar[u] \ar@{->>}[d] & & \\ 
\et & & \etplus \ar[ll] \ar@{->>}[r]^{\theta} &  \OO_{\Cp}/p  & &  \\ } \]
All the rings with tildes ($\tilde{{\quad}}$) also have versions without a tilde: 
one goes from the latter to the former by making Frobenius invertible and completing. 
For example, $\et$ is the completion
of the  perfection of $\e$.

The three rings in the leftmost column (at least their 
tilde-free versions) are related to the theory of 
$(\phi,\Gamma_K)$-modules. 
The three rings in the rightmost column are related 
to $p$-adic Hodge theory. To go from one theory to
the other, one goes from one side to the other through 
all the intermediate rings. The best case is when
one can work in the middle column. For example, from 
top to bottom: semi-stable, crystalline, or finite
height representations. The ring that binds them all is $\btst{}{}$.

\Subsection{List of the rings of power series}
Let us review the different rings of power series which occur in this article; 
let $C[r;1[$
be the annulus$\{ z\in \Cp,\ p^{-1/r} \leq |z|_p < 1 \}$.
We then have:
\[ \begin{array}{||c|c|}
\hline
\mathbf{E}_F^+ & k[[T]] \\
\mathbf{A}_F^+ & \OO_F[[T]] \\
\mathbf{B}_F^+ & F \otimes_{\OO_F} \OO_F[[T]] \\
\hline\hline
\mathbf{E}_F & k((T)) \\
\mathbf{A}_F & \hat{\OO_F[[T]][T^{-1}]} \\
\mathbf{B}_F & F \otimes_{\OO_F} \hat{\OO_F[[T]][T^{-1}]}  \\
\hline\hline
\adag{,r}_F & \text{Laurent series $f(T)$, convergent
on $C[r;1[$, and bounded by $1$} \\
\bdag{,r}_F & \text{Laurent series $f(T)$, convergent 
on $C[r;1[$, and bounded} \\
\hline\hline
\bnrig{,r}{,F} & \text{Laurent series $f(T)$, convergent 
on $C[r;1[$}  \\
\bnst{,r}{,F} & \bnrig{,r}{,F}[\log(T)] \\
\hline\hline
\bhol{,F} & f(T) \in F[[T]],\ 
\text{$f(T)$ converges on the open unit disk $D[0;1[$}  \\
\blog{,F} & \bhol{,F}[\log(T)] \\
\hline
\end{array} \]

\[ \clubsuit\diamondsuit\spadesuit\heartsuit \]

\makelastpage
\end{document}